\documentclass{article}
\usepackage{amsfonts}
\usepackage{amsmath,amsthm}

\usepackage[all]{xy}
\usepackage{graphicx}

\usepackage{amssymb}
\usepackage{latexsym}

\DeclareMathAlphabet\EuFrak{U}{euf}{m}{n}	
\SetMathAlphabet\EuFrak{bold}{U}{euf}{b}{n}	

\newcommand{\ul}{\underline}

\newcommand{\lra}{\leftrightarrow}

\newcommand{\lora}{\longrightarrow}

\newcommand{\hra}{\hookrightarrow}

\newcommand{\ovl}{\overline}

\newcommand{\wa}{\widehat}

\newcommand{\sC}{{\it C*}-}

\newcommand{\bC} {{\mathbb C}}
\newcommand{\bH} {{H}}

\newcommand{\bT} {{\mathbb T}}

\newcommand{\bZ} {{\mathbb Z}}
\newcommand{\bM} {{\mathbb M}}
\newcommand{\bN} {{\mathbb N}}
\newcommand{\bP} {{\mathbb P}}

\newcommand{\ud}{{{\mathbb U}(d)}}
\newcommand{\sud}{{{\mathbb {SU}}(d)}}

\newcommand{\eps}{\varepsilon}

\newcommand{\mA}{\mathcal A}
\newcommand{\mB}{\mathcal B}
\newcommand{\mC}{\mathcal C}
\newcommand{\mD}{\mathcal D}
\newcommand{\mE}{\mathcal E}

\newcommand{\mG}{\mathcal G}

\newcommand{\mI}{\mathcal I}
\newcommand{\mJ}{\mathcal J}

\newcommand{\mL}{\mathcal L}
\newcommand{\mM}{\mathcal M}

\newcommand{\mO}{\mathcal O}
\newcommand{\mP}{\mathcal P}

\newcommand{\mT}{\mathcal T}

\newcommand{\mV}{\mathcal V}

\newcommand{\mW}{\mathcal W}

\newcommand{\tend}{{\bf end}\mA}

\newcommand{\oro}{\mO_\rho}

\newcommand{\autrb}{ {\bf aut}_{\eps_\bullet} \wa \rho_\bullet }

\newcommand{\coe}{{\mO_\mE}}

\newcommand{\rhors}{{\rho^r , \rho^s}}
\newcommand{\hrs}{{ \bH^r , \bH^s }}

\newcommand{\ii}{{\iota,\iota}}

\newcommand{\rr}{{\rho,\rho}}

\newcommand{\rs}{{\rho,\sigma}}
\newcommand{\sr}{{\sigma,\rho}}

\newcommand{\rsp}{{\rho',\sigma'}}

\newcommand{\rt}{{\rho,\tau}}
\newcommand{\tr}{{\tau,\rho}}

\newcommand{\st}{{\sigma,\tau}}
\newcommand{\ts}{{\tau,\sigma}}

\newcommand{\sss}{{\sigma,\sigma}}

\newcommand{\xr}{{\xi,\rho}}
\newcommand{\xs}{{\xi,\sigma}}

\newcommand{\rx}{{\rho,\xi}}
\newcommand{\sx}{{\sigma,\xi}}

\newcommand{\xx}{{\xi,\xi}}

\newcommand{\obc}{{{\bf obj} \ \mC}}
\newcommand{\obt}{{{\bf obj} \ \mT}}

\newcommand{\sym} { {\bf{sym}} }

\newtheorem{thm}{Theorem}[section]
\newtheorem{cor}[thm]{Corollary}
\newtheorem{lem}[thm]{Lemma}
\newtheorem{prop}[thm]{Proposition}

\newtheorem{defn}[thm]{Definition}

\theoremstyle{definition}
\newtheorem{ex}{Example}[section]

\theoremstyle{remark}
\newtheorem{rem}{Remark}[section]

\numberwithin{equation}{section}

\begin{document}

\author{{\sf Ezio Vasselli}
                         \\{\it Dipartimento di Matematica}
                         \\{\it University of Rome "La Sapienza"}
			 \\{\it P.le Aldo Moro, 2 - 00185 Roma - Italy }
                         \\{\sf vasselli@mat.uniroma2.it}}

\title{ Bundles of \sC categories }
\maketitle

\begin{abstract}
We introduce the notions of multiplier \sC category and continuous bundle of \sC categories, as the categorical analogues of the corresponding \sC algebraic notions. Every symmetric tensor \sC category with conjugates is a continuous bundle of \sC categories, with base space the spectrum of the \sC algebra associated with the identity object. We classify tensor \sC categories with fibre the dual of a compact Lie group in terms of suitable principal bundles. This also provides a classification for certain \sC algebra bundles, with fibres fixed-point algebras of $\mO_d$.

\

\end{abstract}

\section{Introduction.}
\label{intro}

\sC categories have been studied in recent years from several points of view, and with several motivations, as duality for (quantum) groups (\cite{DR89,Wor87,DPR01}), quantum field theory (\cite{DR89A,BL04}), $K$-theory (\cite{Mit01}). In the present work, we provide a notion of bundle in the categorical setting, and then show that every \sC category endowed with a tensor structure is a bundle in the above-mentioned sense. This allows us to provide classification results by mean of cohomological methods, and then to approach a duality theory for group bundles and groupoids.

A tensor category is described by a collection of objects $\rho$, $\sigma$, $\ldots$, together with vector spaces denoted by $(\rs)$, $\ldots$, called the {\em spaces of arrows}. Arrows $t \in (\rs)$, $t' \in (\st)$ can be composed to obtain an element $t' \circ t \in (\rt)$, so that in particular every $(\rr)$ is a ring. Moreover, the objects can be multiplied, via the {\em tensor product} $\rs \mapsto$ $\rho \otimes \sigma$, and the existence of an identity, denoted by $\iota$, is postulated in such a way that $\iota \otimes \rho =$ $\rho \otimes \iota =$ $\rho$ for every object $\rho$.
The aim of duality theory is to characterize tensor categories in terms of {\em duals}, i.e. categories with objects vector spaces and arrows linear maps equivariant w.r.t. a (possibly unique) {\em dual object}. In correspondence with the additional structure carried by the given tensor category, the dual object turns out to have different natures; for example, we may get a compact group (\cite{DR89}), an algebraic group (\cite{Del90}), a compact quantum group (\cite{Wor87}), a multiplicative unitary (\cite{DPR01}), just to mention some of the known results.
The above-considered cases are characterized by different properties of the tensor structure: the crucial ones are the symmetry for \cite{DR89,Del90} and the braiding for \cite{Wor87,DPR01}, which describe the commutativity of the tensor product. On the other side, a point which is common to all the cases is the following: the identity object is {\em simple}, in the sense that the ring $(\ii)$ is isomorphic to the field of scalars. 

In the present paper, we proceed with the program started in \cite{Vas} to construct a duality theory for the case in which the given tensor category is symmetric and with non-simple identity. Our setting is the one of {\em tensor \sC categories}, considered in \cite{DR89,Wor87,DPR01}; in particular, this implies that $(\ii)$ is a commutative \sC algebra. The condition that $\iota$ is non-simple is translated as the fact that the spectrum $X^\iota$ of $(\ii)$ does not reduce to a single point. We show that every tensor \sC category has a canonical structure of a bundle, and study the important class of {\em special categories}, which turn out to be bundles with 'fibre' the dual of a compact Lie group. We provide a complete classification for special categories, in terms of a suitable cohomology set. Such a classification admits a natural translation in terms of continuous bundles of \sC dynamical systems, with fibre a fixed-point algebra of the Cuntz algebra $\mO_d$.

Our research program is intended to provide a generalization of the Doplicher-Roberts duality theory (\cite{DR87,DR89,DR90}). Anyway, as we will show in a forthcoming paper, dramatic differences arise in the case with non-simple identity, concerning non-existence and non-unicity of the dual object, which will turn out to be (when existing) a bundle of compact groups, or more generally a groupoid. These phenomena have a cohomological nature, and their root is the classification provided in the present work. As a future application, we will provide a new kind of twisted topological equivariant $K$-functor, related to the gauge-equivariant $K$-theory introduced by V. Nistor and E. Troitsky (\cite{NT04}). The construction of such a $K$-functor is outlined in \cite{Vas06}.

The motivation of \cite{DR89} was the search for the gauge group in the setting of superselection sectors in quantum field theory (\cite{DR90}). Motivated by similar structures arising in low-dimensional quantum field theory and quantum constraints, a duality theory for tensor categories of \sC algebra endomorphisms has been developed in \cite{BL04}: the set $(\ii)$ is identified as the centre of a given 'observable' \sC algebra, and the existence of a subcategory (with simple identity) satisfying the Doplicher-Roberts axioms is postulated, with the result that the reconstructed dual object is a compact group. The tensor \sC categories above considered turn out to be 'trivial bundles' in the sense of the present paper, as remarked in \cite[Ex.5.1]{Vas05A}.

The present work is organized as follows.
In Sec.\ref{sec_mcx}, we introduce the notions of {\em multiplier \sC category} (Prop.\ref{prop_mc}) and $C_0(X)$-{\em category} (Def.\ref{def_cxcat}, $X$ denotes a locally compact, Hausdorff space), which generalize the ones of multiplier \sC algebra and $C_0(X)$-algebra. As an application, the notion of {\em multiplier \sC bimodule} is introduced (Cor.\ref{cor_mbm}). We prove that every $C_0(X)$-category has a canonical bundle structure (Prop.\ref{prop_cx_cat}), and consider the notion of {\em continuous bundle} of \sC categories. 
In Sec.\ref{ssec_tens}, we study tensor \sC categories with $(\ii) \neq \bC$. We prove that every tensor \sC category is a $C(X^\iota)$-category in a canonical way, where $C(X^\iota)$ is a suitable \sC subalgebra of $(\ii)$ (Prop.\ref{prop_tens_cx}). In particular, a symmetric tensor \sC category with (twisted) conjugates is a {\em continuous} bundle, having as fibres duals of compact groups (Thm.\ref{thm_tcf}, Rem.\ref{rem_fib_dual}). The above-quoted theorem has an analogue in the setting of $2$-\sC categories with conjugates (\cite[\S 3]{Zit05}); our proof is different, and is based on the notion of {\em (twisted) special object} (see (\ref{eq_tso})). Once that a bundle structure is considered for symmetric tensor \sC categories, it is natural to formulate a coherent notion of {\em local triviality} (Def.\ref{def_slc}). A {\em special category} is a locally trivial, symmetric tensor \sC category such that the fibre $\mT_\bullet$ is isomorphic (up to direct sums) to the dual of a compact Lie group $G$ (Def.\ref{def_spec_cat}). Special categories with fibre $\mT_\bullet$ are completely classified by the cohomology set $H^1(X^\iota,QG)$ (Thm.\ref{thm_amen}), where $QG$ is a suitable compact Lie group associated with $G$ (see (\ref{def_p})). Finally, we consider a canonical $G$-action on $\mO_d$, and give a classification for continuous bundles with fibre the fixed-point algebra $\mO_G$ (Sec.\ref{ss_ogb}).

\subsection{Notation and Keywords.}

Let $X$ be a locally compact Hausdorff space. We denote by: $C_0(X)$, the \sC algebra of ($\bC$-valued) continuous functions on $X$ vanishing at infinity; $C_b(X)$, the \sC algebra of bounded, continuous functions on $X$; $C(X)$, for $X$ compact, the \sC algebra of continuous functions on $X$; $C_U(X)$, for $U \subseteq X$ open, the ideal in $C_0(X)$ of functions vanishing on $X-U$; $C_x(X) :=$ $C_{X - \left\{ x \right\} } (X)$, $x \in X$.
If $\left\{ X_i \right\}$ is a cover of $X$, then we denote $X_{ij} := X_i \cap X_j$, $X_{ijk} := X_i \cap X_j \cap X_k$.

If $\mA$ is a \sC algebra, we denote by ${\bf aut} \mA$ (resp. $\tend$) the set of automorphisms (resp. endomorphisms) of $\mA$, endowed with the topology of pointwise convegence. A pair $( \mA , \rho)$, where $\rho \in \tend$, is called \sC dynamical system. If $(\mA , \rho)$, $(\mA' , \rho')$ are \sC dynamical systems, a \sC algebra morphism $\alpha : \mA \to \mA'$ such that $\alpha \circ \rho = \rho' \circ \alpha$ is denoted by $\alpha : ( \mA , \rho ) \to ( \mA' , \rho' )$.

Finally, for every $d \in \bN$ we denote by $\ud$ the unitary group, and by $\sud$ the special unitary group.

\subsubsection{$C_0(X)$-algebras.}
\label{ssec_cxalg}
Let $X$ be a locally compact Hausdorff space. A $C_0(X)${\em -algebra} is a \sC algebra $\mA$ endowed with a non-degenerate morphism from $C_0(X)$ into the centre of the multiplier algebra $M(\mA)$ (\cite[\S 2]{Kas88}). It is customary to identify elements of $C_0(X)$ with their images in $M(\mA)$. We call $C_0(X)$-morphism a \sC algebra morphism commuting with the above $C_0(X)$-action; in particular, we denote by ${\bf aut}_X \mA$ (resp. ${\bf end}_X \mA$) the set of $C_0(X)$-automorphisms (resp. $C_0(X)$-endomorphisms) of $\mA$. It can be proved that $C_0(X)$-algebras correspond to upper-semicontinuous bundles of \sC algebras with base space $X$ (\cite{Nil96}); in particular, every continuous bundle of \sC algebras (\cite[Chp.10]{Dix},\cite{KW95}) is a $C_0(X)$-algebra. The fibres of $\mA$ are defined as the quotients $\mA_x := \mA / (C_x (X) \mA)$, $x \in X$. We denote by $\otimes_X$ the minimal $C_0(X)$-algebra tensor product \cite[\S 2]{Kas88}. If $\mA$, $\mB$ are $C_0(X)$-algebras, a $C_0(X)$-{\em Hilbert $\mA$-$\mB$-bimodule} is a Hilbert $\mA$-$\mB$-bimodule $\mM$ such that $af\psi b = a \psi f b$, $f \in C_0(X)$, $a \in \mA$, $b \in \mB$, $\psi \in \mM$. 

\subsubsection{bi-Hilbertian bimodules.}
\label{ss_bh}
The following notion appeared in \cite{KPW04}. Let $\mA$, $\mB$ be \sC algebras, $\mM$ a Hilbert $\mA$-$\mB$-bimodule such that the \sC algebra $K(\mM)$ of compact, right $\mB$-module operators is contained in $\mA$. Then, besides the usual $\mB$-valued scalar product $\left \langle \cdot , \cdot \right \rangle_\mB$, there exists an $\mA$-valued scalar product $\left \langle \psi , \psi' \right \rangle_\mA := \theta_{\psi , \psi'}$, where $\psi , \psi' \in \mM$, and $\theta_{\psi,\psi'} \in K(\mM) \subseteq \mA$, $\theta_{\psi , \psi'} \varphi := \psi \left \langle \psi' , \varphi \right \rangle_\mB$, $\varphi \in \mM$. We say that $\mM$ is a {\em bi-Hilbertian bimodule} if $\left\| \psi \right\|^2 = \left\| \theta_{\psi , \psi} \right\|$, $\psi \in \mM$, and use the notation $_\mA \mM_\mB$. A morphism of bi-Hilbertian bimodules is a bounded linear map $\phi : \ _\mA\mM_\mB \to \ _{\mA'}\mM'_{\mB'}$, with \sC algebra morphisms $\phi_l : \mA \to \mA'$, $\phi_r : \mB \to \mB'$ such that $\phi (a \psi b) = \phi_l (a) \phi (\psi) \phi_r (b)$, $\phi_r (\left \langle  \psi , \psi' \right \rangle_\mB) = \left \langle \phi (\psi) , \phi (\psi') \right \rangle_{\mB'}$, $\phi_l ( \left \langle  \psi , \psi' \right \rangle_\mA ) = \left \langle \phi (\psi) , \phi (\psi') \right \rangle_{\mA'}$, $\psi , \psi' \in \mM$, $a \in \mA$, $b \in \mB$.

\subsubsection{Bundles.}
\label{ss_bun}
For standard notions about {\em vector bundles}, we refer to the classics \cite{Ati,Kar,Seg68}. In the present paper, we will assume that every {\em vector bundle is endowed with a Hermitian structure} (see \cite[I.8]{Kar}). If $X$ is a locally compact Hausdorff space, we denote by ${\bf vect}(X)$ the category having as objects vector bundles with base space $X$, and arrows vector bundle morphisms. In the present work, we will also deal with {\em Banach bundles} in the sense of \cite{Dup74} (i.e., bundles of Banach spaces that are not necessarily locally trivial). About the related notion of continuous field of Banach spaces, we refer to \cite[Chp.10]{Dix}.

\section{Multipliers, $C_0(X)$-categories.}
\label{sec_mcx}

A category $\mC$ not necessarily endowed with identity arrows is said to be a {\em \sC category} if every set of arrows $(\rs)$, $\rs \in \obc$, is a Banach space such that the composition of arrows defines linear maps satisfying $\left\| t' \circ t \right\| \leq \left\| t' \right\| \left\| t \right\|$, $t \in (\rs)$, $t' \in (\st)$; moreover, an antilinear isometric cofunctor $* : \mC \to \mC$ is assigned, in such a way that $\left\| t^* \circ t \right\| = \left\| t \right\|^2$. 
Our \sC category $\mC$ is said to be {\em unital} if every $\rho \in \obc$ admits an identity arrow $1_\rho \in (\rr)$ such that $t \circ 1_\rho = t$, $t \in (\rs)$. Functors between \sC categories preserving the above structure are said {\em \sC functors}; in particular, we will use the term {\em \sC monofunctor}, resp. \sC {\em autofunctor}, resp. {\em \sC epifunctor} in the case in which the given \sC functor induces injective, resp. one-to-one, resp. surjective Banach space maps on the spaces of arrows. If $\mD$ is a \sC subcategory of $\mC$, then we use the notation $\mD \subseteq \mC$; in such a case, every space $(\rs)_\mD$ of arrows in $\mD$ is a Banach subspace of $(\rs)$. 
Basic references for \sC categories are \cite{DR89,KPW04}; non-unital \sC categories have been considered in \cite{Mit01}. 

For the notions of {\em subobjects} and {\em direct sums}, we refer to \cite[\S 1]{DR89}.
The {\em closure for subojects} of a \sC category $\mC$ is the \sC category $\mC_s$ with objects $E = E^* = E^2 \in$ $(\rr)$, $\rho \in \obc$, and arrows $(E,F) :=$ $\left\{ t \in (\rs) : t = F \circ t = t \circ E \right\}$; by construction, $\mC_s$ has subobjects. The {\em additive completition} of $\mC$ is the \sC category $\mC_+$ with objects $n$-ples $\ul \rho := ( \rho_1 , \ldots , \rho_n )$, $n \in \bN$, and arrows spaces of matrices $(\ul \rho , \ul \sigma) :=$ $\left\{ (t_{ij}) : t_{ij} \in (\rho_j , \sigma_i ) \right\}$; by construction, $\mC_+$ has direct sums (see \cite[Def.2.12]{Mit01}).

A (unital) \sC category with a single object is a (unital) \sC algebra. A \sC category with two objects is a bi-Hilbertian \sC bimodule, together with the conjugate bimodule. 

\subsection{Multipliers.}
\label{categories}

Let $\mC$ be a \sC category. Then, every $(\rr)$, $\rho \in \obc$, is a \sC algebra, and every $(\rs)$ can be endowed with the following structure of a bi-Hilbertian $(\sss)$-$(\rr)$-bimodule:
\[
\left\{
\begin{array}{ll}
b , t \mapsto b \circ t    \ \  , \ \ 
t , a \mapsto t \circ a     \\
\left \langle t , t' \right \rangle_\rho := t^* \circ t' \in (\rr)  \ , \ \
\left \langle t' , t \right \rangle_\sigma := t' \circ t^* \in (\sss) \ ,
\end{array}
\right.
\]
\noindent $a \in (\rr)$, $b \in (\sss)$, $t,t' \in (\rs)$; in fact, $\left\| t  \right\|^2 = \left\| \right. \left \langle t , t \right \rangle_\rho \left. \right\| = \left\| \left \langle t , t \right \rangle_\sigma  \right\|$. Let $\rs , \xi \in \obc$; according to the above-defined structure, we say that a bounded linear map $T : (\xr) \to (\xs)$ is a {\em right $(\xx)$-module operator} if $T (t \circ a) = T(t) \circ a$, $t \in (\xr)$, $a \in (\xx)$. In the same way, a bounded linear map $T' : (\rx) \to (\sx)$ is a {\em left $(\xx)$-module operator} if $T' (a \circ t) = a \circ T'(t)$, $t \in (\rx)$, $a \in (\xx)$.
A {\em multiplier} from $\rho$ to $\sigma$ is a pair $( T^l , T^r )$, where $T^l : (\rr) \to (\rs)$ is a right $(\rr)$-module operator and $T^r : (\sss) \to (\rs)$ is a left $(\sss)$-module operator, such that the relation
\begin{equation}
\label{eq_dm}
b \circ T^l (a) = T^r (b) \circ a
\end{equation}
\noindent is satisfied for every $a \in (\rr)$, $b \in (\sss)$. We denote by $M(\rs)$ the set of multipliers from $\rho$ into $\sigma$. Note that (\ref{eq_dm}) implies $\left\| T^l  \right\| = \left\| T^r  \right\|$. In order to economize in notation, we define
\[
T a := T^l (a)  \ , \ \ b T := T^r (b) \ , \ \ T \in M(\rs) \ ,
\]
\noindent so that (\ref{eq_dm}) can be regarded as an associativity constraint:
\begin{equation}
\label{eq_ass}
b \circ Ta = bT \circ a \ \in (\rs)  \ .
\end{equation}

\begin{ex}
\label{ex_arr}
Let $t \in (\rs)$. We define $M_t a := t \circ a$, $b M_t := b \circ t$, $a \in (\rr)$, $b \in (\sss)$. It is clear that $M_t \in M(\rs)$.
For every $\rho \in \obc$, we also define $1_\rho a = a 1_\rho := a$, $a \in (\rr)$. It is clear that $1_\rho \in M(\rr)$; we call $1_\rho$ the {\bf identity multiplier} of $\rho$.
\end{ex}

An involution is defined on multipliers, by 
\begin{equation}
\label{def_*}
\begin{array}{ll}
M(\rs) \to M(\sr) \\
T \mapsto T^* \ : \ 
T^* b := (b^*T)^*
\ , \ 
a T^* := (T a^*)^* \ .
\end{array}
\end{equation}
\noindent In the case $\rho = \sigma$, a composition law can be naturally defined: if $A , A' \in M(\rr)$, then the maps
\[
(\rr) \ni a \mapsto A'(A a) \ , \ \ (\rr) \ni a \mapsto (a A)A' \ ,
\]
\noindent define a multiplier $A' \circ A \in M(\rr)$. It is now clear that every $A \in M(\rr)$ defines a multiplier of $(\rr)$ in the usual \sC algebra sense, so that $M(\rr)$, endowed with the above *-algebra structure, coincides with the multiplier algebra of $(\rr)$, with identity $1_\rho$.

\begin{lem}
\label{lem_m}
Let $\mC$ be a \sC category. For every $\rs , \tau \in \obc$, there are maps
\begin{equation}
\label{eq_cl}
\begin{array}{ll}
M (\st) \times (\rs) \to (\rt)
\ , \ \ 
(S,t) \mapsto St := \lim_\lambda (S e^\sigma_\lambda) \circ t
\\
(\st) \times M (\rs) \to (\rt)
\ , \ \ 
(t',T) \mapsto t'T := \lim_\lambda t' \circ (e^\sigma_\lambda T)
\end{array}
\end{equation}
\noindent where $(e^\sigma_\lambda)_\lambda \subset (\sss)$ is an approximate unit. The above maps naturally extend the composition of arrows in $\mC$.
\end{lem}

\begin{proof}
By \cite[Prop.2.16]{KPW04}, every $(\rs)$ is non-degenerate as a Hilbert bimodule w.r.t. the left $(\sss)$-action and right $(\rr)$-action. Thus, every $t \in (\rs)$ admits factorizations $t = b \circ t_2 = t_1 \circ a$, $t_1 , t_2 \in (\rs)$, $a \in (\rr)$, $b \in (\sss)$ (see \cite[Prop.1.8]{Bla96}). Let $S \in M(\st)$. We consider the net $\left\{ (S e^\sigma_\lambda) \circ t \right\}_\lambda$, and estimate
\[
\begin{array}{ll}
\left\|  (S e^\sigma_\lambda) \circ t - (Sb) \circ t_2           \right\|     & =
\left\|  S ( e^\sigma_\lambda \circ b - b ) \circ t_2            \right\|  \\ & \leq
\left\|  S ( e^\sigma_\lambda \circ b - b ) \right\| \left\| t_2 \right\|  \\ & \leq

\left\| S \right\| 
\left\| e^\sigma_\lambda \circ b - b  \right\|
\left\| t_2 \right\| \ .
\end{array}
\]
\noindent Since by definition of approximate units $\left\| e^\sigma_\lambda \circ b - b  \right\| \stackrel{\lambda}{\to} 0$, we conclude that $\left\{ (S e^\sigma_\lambda) \circ t \right\}_\lambda$ converges to $(Sb) \circ t_2 \in (\rs)$. Note that the limit is unique, and that it does not depend on the choice of the approximate unit. The same argument applies for the net $\left\{ t' \circ (e^\sigma_\lambda T) \right\}$, and this proves the lemma. 
\end{proof}

We can now introduce the composition law
\begin{equation}
\label{def_cl}
\begin{array}{ll}
M (\st) \times M(\rs) \to M(\rt) \\
S , T \mapsto S \circ T \ : \
(S \circ T) \ a := S ( Ta ) \ , \ 
c \ (S \circ T) := (cS) T \ ,

\end{array}
\end{equation}
\noindent $a \in (\rr)$, $c \in (\tau , \tau)$, $T \in M(\rs)$, $T' \in M(\st)$. Note in fact that $Ta \in (\rs)$, $cS \in (\tau , \sigma)$, thus we can apply the previous lemma and perform the compositions $S ( Ta )$, $(cS) T$.

A \sC subcategory $\mI \subseteq \mC$, with arrows $(\rs)_\mI \subseteq (\rs)$, $\rs \in {\bf obj} \mI = {\bf obj} \mC$, is said to be an {\em ideal} \ if $t \circ t' \in (\ts)_\mI$, $t'' \circ t \in (\rx)_\mI$ for every $t \in (\rs)_\mI$, $t' \in (\tr)$, $t'' \in (\sx)$. If $\mI \subseteq \mC$ is an ideal, we write $\mI \lhd \mC$. In particular, $\mI$ is said to be an {\em essential ideal} if every ideal of $\mC$ has nontrivial intersection with $\mI$ (i.e., $(\rs)_\mI \cap (\rs)_\mJ \neq \left\{ 0 \right\}$ for some $\rs \in \obc$, $\mJ \lhd \mC$).

\begin{prop}
\label{prop_mc}
Let $\mC$ be a \sC category. The category $M(\mC)$ having the same objects as $\mC$, and arrows $M(\rs)$, $\rs \in \obc$, is a unital \sC category. Moreover, the maps $$\left\{ (\rs) \ni t \mapsto M_t \in M(\rs) \right\}$$ defined in Ex.\ref{ex_arr} induce a \sC monofunctor $I : \mC \hra M(\mC)$, in such a way that $I(\mC)$ is an essential ideal of $M(\mC)$. If $\mC$ is unital, then $I(\mC) = M(\mC)$. If $\mC$ is an essential ideal of a unital \sC category $\mC'$, then there exists a \sC monofunctor $I' : \mC' \hra M(\mC)$ extending $I$.
\end{prop}

\begin{proof}
It follows from Lemma \ref{lem_m} that $M(\mC)$ is a \sC category, with units defined as by Ex.\ref{ex_arr}. It is clear that $I$ is a \sC functor; by well-known properties of multiplier algebras the maps $I : (\rr) \to M(\rr)$ are isometric, thus if $t \in (\rs)$, then $\left\| t \right\|^2$ $=$ $\left\| t^* \circ t \right\|$ $=$ $\left\| M (t^* \circ t) \right\|$ $=$ $\left\| M(t) \right\|^2$. This implies that $I$ is a \sC monofunctor. Let now $\mC'$ be a unital \sC category as above; we denote by $(\rs)'$ the spaces of arrows of $\mC'$ (by definition of essential ideal, we identify $\obc$ with $\obc'$). For every $w \in (\rs)$ we define a multiplier $M_w \in M (\rs)$,
\[
\left\{
\begin{array}{ll}
M_w a := w \circ a  \ , \ \ a \in (\rr)
\\
b M_w := b \circ w  \ , \ \ b \in (\sss)
\end{array}
\right.
\]
\noindent (note in fact that $w \circ a$, $b \circ w$ $\in$ $(\rs)$). In this way, we obtain a \sC functor $I' : \mC' \to M(\mC)$, $w \mapsto M_w$, which reduces to $I$ for arrows in $\mC$. In the case $\rho = \sigma$, we find that $I'$ induces \sC algebra morphisms $I'_\rho : ( \rr )' \to M(\rr)$; such morphisms are injective by \cite[Prop.3.12.8]{Ped}. Thus, if $w \in (\rs)'$, then $\left\| w \right\|^2$ $=$ $\left\| w^* \circ w \right\|$ $=$ $\left\| M_{w^* \circ w} \right\|$ $=$ $\left\| M_w \right\|^2$. 
\end{proof}

In the case of \sC categories with two objects, we obtain a construction for bi-Hilbertian bimodules:

\begin{cor}
\label{cor_mbm}
Let $\mA$, $\mB$ be \sC algebras, $_\mA \mM_\mB$ a bi-Hilbertian bimodule. Then, there exists a universal bi-Hilbertian $M(\mA)$-$M(\mB)$-bimodule $M(\mM)$ with a monomorphism $I : \ _\mA \mM_\mB \to \ _{M(\mA)} M(\mM)_{M(\mB)}$, and satisfying the universal property w.r.t. bi-Hilbertian bimodules containing $\mM$ as essential ideal.
\end{cor}

\begin{ex}
\label{ex_mvb}
Let $X$ be a locally compact, paracompact space, $d \in \bN$, and $\mE \to X$ a rank $d$ vector bundle. We denote by $\Gamma_0 \mE$ the $C_0(X)$-module of continuous sections of $\mE$ vanishing at infinity, and by $\Gamma_0 (\mE , \mE)$ the \sC algebra of compact operators of $\Gamma_0 \mE$. Then, $\Gamma_0 (\mE , \mE)$ is a continuous bundle of \sC algebras over $X$, with fibre the matrix algebra $\bM_d$; moreover, $\Gamma_0 \mE$ naturally becomes a Hilbert $\Gamma_0(\mE,\mE)$-$C_0(X)$-bimodule. We consider the \sC category $\mC$ with objects $\rs$, and arrows $(\rr) := C_0(X)$, $(\sss) := \Gamma_0 (\mE , \mE )$, $(\rs) := \Gamma_0 (\mE)$, $(\sr) := \Gamma_0^*\mE$ (i.e., the conjugate bimodule of $\Gamma_0(\mE)$). Let $\Gamma_b (\mE , \mE)$ denote the \sC algebra of bounded, continuous sections of the \sC bundle associated with $\Gamma_0 (\mE , \mE)$, and $\Gamma_b \mE$ the module of bounded, continuous sections of $\mE$. Then, $M(\mC)$ has arrows $M(\rr) = C_b(X)$, $M(\sss) := \Gamma_b (\mE , \mE)$, $M(\rs) = \Gamma_b \mE$ (see \cite[Thm.3.3]{APT73}).
\end{ex}

A different approach to multiplier \sC categories can be found in \cite{Kan01}. Since the universal property Prop.\ref{prop_mc} is satisfied by the \sC categories introduced in the above-cited reference, the two constructions provide the same result.

\subsection{$C_0(X)$-categories.}
\label{sec_cx_cat}
Let $\mC$ be a \sC category. For every $\rho \in \obc$, we denote by $ZM(\rr)$ the centre of the \sC algebra $M(\rr)$ of multipliers of $(\rr)$. 

\begin{defn}
\label{def_cxcat}
Let $\mC$ be a \sC category, $X$ a locally compact Hausdorff space. $\mC$ is said to be a {\bf $C_0(X)$-category} if for every $\rho \in \obc$ there is given a non-degenerate morphism $i_\rho : C_0(X) \to M(\rr)$, such that for every $\sigma \in \obc$, $t \in (\rs)$, $f \in C_0(X)$, the following equality holds:
\begin{equation}
\label{eq_cxcat}
t \ i_\rho(f) = i_\sigma (f) \ t \ .
\end{equation}
\end{defn}

The set $\left\{  i_\rho \right\}_{\rho \in \obc}$ is called $C_0(X)$-{\em action}. It follows from (\ref{eq_cxcat}) that $i_\rho (f) \in ZM (\rr)$ for every $\rho \in \obc$, thus $(\rr)$ is a $C_0(X)$-algebra. By definition, every $(\rs)$ is a $C_0(X)$-Hilbert $(\sss)$-$(\rr)$-bimodule in the sense of Sec.\ref{ssec_cxalg}. A \sC functor $\eta : \mC \to \mC'$ between $C_0(X)$-categories is said to be a $C_0(X)$-functor if $\eta ( t \ i_\rho (f) ) = \eta (t) \ i'_{\eta(\rho)} (f)$ for every $t \in (\rs)$, $\rs \in \obc$, $f \in C_0(X)$. In the sequel, we will drop the symbol $i_\rho$, so that (\ref{eq_cxcat}) simply becomes $tf = ft$. 

\subsubsection{The bundle structure.}
Let $X$ be a locally compact Hausdorff space, $\mC$ a $C_0(X)$-category, $U \subseteq X$ an open set. The {\em restriction of $\mC$ over $U$} is the \sC subcategory $\mC_U$ having the same objects as $\mC$, and arrows the Banach spaces $(\rs)_U$ $:=$ $C_U(X) (\rs)$ $:=$ ${\mathrm{span}}$ $\left\{ f t , f \in C_U (X) , t \in (\rs) \right\}$; it is clear that $\mC_U \lhd \ \mC$ is an ideal, and that $\mC_U$ is a $C_0(U)$-category. 
Let now $W \subseteq X$ be a closed set. 
The {\em restriction of $\mC$ over $W$} is defined as the \sC category $\mC_W$ having the same objects as $\mC$, and arrows the quotients $(\rs)_W :=$ $(\rs) / [C_{X-W}(X)(\rs)]$. It is easily checked that $\mC_W$ is a \sC category, in fact composition of arrows and involution factorize through elements of $C_{X-W}(X) (\rs)$, $\rs \in \obc$. 
By definition, there is an exact sequence of \sC functors
\begin{equation}
\label{eq_es_restr}
{\bf 0} \lora
\mC_{X-W} 
\stackrel{i_{X-W}}{\lora} 
\mC 
\stackrel{\eta_W}{\lora} 
\mC_W 
\ .
\end{equation}
In the case in which $C_b(W)$ is a $C_0(X)$-algebra (for example, in the case of $X$ being normal), then $\mC_W$ may be described as the \sC category with arrows $ (\rs) \otimes_X C_b(W)$, where $\otimes_X$ is the tensor product with coefficients in $C_0(X)$ defined as in Sec.\ref{ssec_cxalg}. Moreover, $\mC_W$ has an obvious structure of $C(\beta W)$-category, where $\beta W$ denote the Stone-Cech compactification of $W$.

\begin{defn}
\label{def_fibre}
Let $\mC$ be a $C_0(X)$-category, $x \in X$. The {\bf fibre} of $\mC$ over $x$ is defined as the restriction $\mC_x := \mC_{\left\{ x \right\}}$. The restriction \sC epifunctor $\pi_x : \mC \to \mC_x$ is called the {\bf fibre functor}. We denote by $(\rs)_x$ the spaces of arrows of $\mC_x$, $x \in X$.
\end{defn}

\begin{rem}
\label{rem_str_fibre}
Let $\mC$ be a \sC category. Two objects $\rs \in \obc$ are said to be unitarily equivalent if there exists $u \in (\rs)$ such that $u \circ u^* = 1_\sigma$, $u^* \circ u = 1_\rho$. We denote by ${\bf obj_u} \mC$ the set of unitary equivalence classes of objects of $\mC$.
Let now $\mC$ be a $C_0(X)$-category. By definition, the restrictions of $\mC$ over open (closed) subsets of $X$ have the same objects as $\mC$. In particular, this is true for the fibres $\mC_x$, $x \in X$.
Now, it is clear that the fibre functors map unitarily equivalent objects into unitarily equivalent objects. On the other side, $\rs$ may be unitarily equivalent as objects of every $\mC_x$, $x \in X$, but this does not imply that $\rs$ are unitarily equivalent in $\mC$.
Thus, the maps ${\bf obj_u} \mC \to {\bf obj_u} \mC_x$ induced by the fibre functors are always surjective, but not injective in general.
\end{rem}

Let $F : \mC \to \mC'$ be a $C_0(X)$-functor. Then, for every $x \in X$ there exists a \sC functor $F_x : \mC_x \to \mC'_x$ such that 
\begin{equation}
\label{eq_ff}
F_x \circ \pi_x = \pi'_x \circ F
\end{equation}
\noindent where $\pi'_x : \mC' \to \mC'_x$ are the fibre functors defined on $\mC'$. In fact, if $t,v \in (\rs)$ and $\pi_x (t-v) = 0$, then there exists a factorization $t-v = f w$, $f \in C_x(X)$, $w \in (\rs)$, so that $\pi'_x \circ F (t-v) = f(x) \pi'_x \circ F (w) = 0$, and the l.h.s. of (\ref{eq_ff}) is well-defined. We have the following categorical analogue of \cite[Thm.2.3]{Nil96},\cite[Cor.1.12,Prop.2.8]{Bla96}.
\begin{prop}
\label{prop_cx_cat}
Let $X$ be a locally compact Hausdorff space, $\mC$ a $C_0(X)$-category, and $\left\{ \pi_x : \mC \to \mC_x \right\}_{x \in X}$ the set of fibre functors of $\mC$. For every $\rs \in \obc$, $t \in (\rs)$, the norm function $n_t(x) := \left\| \pi_x (t) \right\|$, $x \in X$, is upper semicontinuous and vanishing at infinity; moreover, $\left\| t \right\| = \sup_{x \in X} \left\| \pi_x(t) \right\|$.
\end{prop}

\begin{proof}
Let $t \in (\rs)$, $\rs \in \obc$, so that $t^* \circ t \in (\rr)$. Since $(\rr)$ is a $C_0(X)$-algebra w.r.t. the $C_0(X)$-action of Def.\ref{def_cxcat}, it follows from \cite[Thm.2.3]{Nil96} that the map $n_t (x) = \left\|  \pi_x (t^* \circ t)  \right\|^{1/2}$, $x \in X$, is upper-semicontinuous; moreover, $\left\| t \right\|^2 =$ $\left\| t^* \circ t  \right\| = $ $\sup_x \left\|  \pi_x (t^* \circ t)  \right\|$. 
\end{proof}

Let $\mC$ be a $C_0(X)$-category. If the norm function $n_t$ is continuous for every arrow $t \in (\rs)$, $\rs \in \obc$, then we say that $\mC$ is a {\em continuous bundle of} \sC categories. In such a case, it follows from the above definition that the spaces of arrows of $\mC$ are continuous fields of Banach spaces; in particular, $(\rr)$ is a continuous bundle of \sC algebras for every $\rho \in \obc$. The proof of the following lemma is trivial, thus it will be omitted.

\begin{lem}
\label{lem_ccx}
Let $X$ be a locally compact Hausdorff space, $\mC$ a continuous bundle of \sC categories over $X$. Then, the closure for subobjects $\mC_s$ and the additive completition $\mC_+$ are continuous bundles of \sC categories.
\end{lem}

\begin{rem}
\label{rem_rest}
Let $X$ be a locally compact Hausdorff space, $\mC$ a continuous bundle of \sC categories over $X$, $\rs \in \obc$. In the definition of continuous field by Dixmier-Douady, the vanishing-at-infinity property is not assumed for norm functions; as a consequence of this fact, if $U \subset X$ is open and $W := \ovl U$ is the closure, then the restriction $\eta_W (\rs)$ defines a unique continuous field over $U$, i.e. the restriction $\eta_W (\rs) |_U$ in the sense of \cite[10.1.7]{Dix}. In order to be concise, we define $\dot{\eta}_W (\rs) := \eta_W (\rs) |_U$.
\end{rem}

%
%
%
%
%

\subsubsection{Local triviality and bundle operations.}
\label{ss_lt_bo}
Let $\mC_\bullet$ be a \sC category; we define the {\em constant bundle} $X \mC_\bullet$ as the $C_0(X)$-category having the same objects as $\mC_\bullet$, and arrows the spaces $(\rs)^X$ of continuous maps from $X$ into $(\rs)$ vanishing at infinity, $\rs \in \obc_\bullet$ (the structure of $C_0(X)$-category is defined in the obvious way). In the case in which $X$ is compact, we consider the spaces of continuous maps from $X$ into $(\rs)$; if $t \in (\rs)$, then with an abuse of notation we will denote by $t \in (\rs)^X$ the corresponding constant map. 

Let $\mC$ be a $C_0(X)$-category, $\mC_\bullet$ a \sC category. We say that $\mC$ is {\em locally trivial with fibre} $\mC_\bullet$ if for every $x \in X$ there exists a closed neighborhood $W \ni x$ with a $C_b(W)$-isomorphism $\alpha_W : \mC_W \to W \mC_\bullet$, in such a way that the induced map $\alpha_W : \obc \to {\bf obj} \ \mC_\bullet$ does not depend on the choice of $W$. The functors $\alpha_W$ are called {\em local charts}.

The above definition implies that $\rho_\bullet := \alpha_W (\rho)$, $\rho \in \obc$, does not depend on $W$; in this way, every space of arrows $(\rs)$ is a locally trivial continuous field of Banach spaces with fibre $( \rho_\bullet , \sigma_\bullet )$, trivialized over subsets which do not depend on the choice of $\rs$. It is clear that we may give the analogous definition by using open neighborhoods, anyway for our purposes it will be convenient to use closed neighborhoods, as the corresponding local charts map {\em unital} \sC categories into {\em unital} \sC categories.

\begin{ex}
Let $\mA$ be a $C_0(X)$-algebra. The \sC category with objects the projections of $\mA$, and arrows $(E,F) := \left\{ t \in \mA \ : \ t = Ft = tE \right\}$ is a $C_0(X)$-category.
\end{ex}

\begin{ex}
Let $\mC$ be a unital \sC category, $X$ a compact Hausdorff space. The extension \sC category $\mC^X$ is defined as the closure for subobjects of the constant bundle $X \mC$. It turns out that $\mC$ is a continuous bundle of \sC categories w.r.t. the $C(X)$-structure induced by $X \mC$ (see \cite{Vas}).
\end{ex}

\begin{ex}
\label{ex_vx}
Let $X$ be a locally compact Hausdorff space. Then, ${\bf vect}(X)$ (Sec.\ref{ss_bun}) is a continuous bundle of \sC categories with fibre the category ${\bf hilb}$ of Hilbert spaces. If $X$ is locally contractible, then ${\bf vect}(X)$ is locally trivial; in fact, if $W \subseteq X$ is contractible then every vector bundle $\mE \to X$ can be trivialized on $W$. If $X$ is not locally contractible, then in general ${\bf vect}(X)$ may be not locally trivial; in fact, it could be not possible to trivialize all the elements of ${\bf vect}(X)$ over the same closed subset.
\end{ex}

Let $X$ be a locally compact, paracompact Hausdorff space. We assume that there exists a locally finite, open cover $\left\{ X_i \right\}$ and, for every index $i$, a continuous bundle of \sC categories $\mC_i$ over the closure $\ovl X_i$, with $C_0(X_{ij})$-isofunctors 
\[
\alpha_{ij} : \eta_{j,X_{ij}} (\mC_j ) \to \eta_{i,X_{ij}} (\mC_i)  
\ \ : \ \
\alpha_{ij} \circ \alpha_{jk} = \alpha_{ik}
\ .
\]
\noindent In the previous expression, $\eta_{i,X_{ij}}$ denotes the restriction epifunctor of $\mC_i$ over $\ovl X_{ij}$. By hypothesis, the sets $\left\{ {\bf obj} \ \mC_i \right\}_i$ are in one-to-one correspondence. We want to define a category $\mC^b$, having the same objects as a fixed $\mC_k$, and arrows the {\em glueings} of the continuous fields of Banach spaces $\left\{ ( \alpha_{ik}(\rho) , \alpha_{ik}(\sigma) )  \right\}_i$, $\rs \in {\bf obj} \ \mC_k$, in the sense of \cite[10.1.13]{Dix} (note that $\alpha_{ik}(\rho) , \alpha_{ik}(\sigma) \in$ ${\bf obj} \ \mC_i$, thus $( \alpha_{ik}(\rho) , \alpha_{ik}(\sigma) )$ is a continuous field over $\ovl X_i$). Since in the above reference {\em open} covers are used, for every index $i$ we consider the restriction $( \alpha_{ik}(\rho) , \alpha_{ik}(\sigma) )  |_{X_i}$, and define the spaces of arrows $(\rs)^b$ of $\mC^b$ as the glueings of  
$\left\{ \ ( \alpha_{ik}(\rho) , \alpha_{ik}(\sigma) )  |_{X_i} \ \right\}_i$
w.r.t. the isomorphisms
\[
\alpha_{ij} \ : \
\dot{\eta}_{j,X_{ij}} (  \alpha_{jk}(\rho) , \alpha_{jk}(\sigma)  )
\ \to \
\dot{\eta}_{i,X_{ij}} (  \alpha_{ik}(\rho) , \alpha_{ik}(\sigma)  )
\]
\noindent (recall Rem.\ref{rem_rest}). Let $\rs \in \obc^b$; by definition, the elements $t \in (\rs)^b$ are in one-to-one correspondence with families 
\begin{equation}
\label{def_glue}
(t_i)_i \in \prod_i ( \ \alpha_{ik} (\rho) , \alpha_{ik} (\sigma) \ )  \ : \ 
\alpha_{ij} \circ \eta_{j,X_{ij}} (t_j)  = \eta_{i,X_{ij}} (t_i)
\ .
\end{equation}
\noindent We define an involution $(t_i)^* := (t_i^*)$ and a composition $(t'_i) \circ (t_i) := (t'_i \circ t_i)$, in such a way that $\mC^b$ is a \sC category. 
Now, our definition does not ensure that $\mC^b$ is a $C_0(X)$-category, in fact the $C_0(X)$-action on a continuous field of Banach spaces may be degenerate in the case of $X$ being not compact. We define a \sC category $\mC$ having the same objects as $\mC^b$, and arrows the spaces $(\rs) :=$ $C_0(X) (\rs)^b$. It is clear that $\mC$ is an ideal of $\mC^b$. We call $\mC$ the {\em glueing} of the categories $\mC_i$. By construction, $\mC$ is a continuous bundle on $X$, satisfying the following universal property: for every index $i$ there exists a \sC isofunctor $\pi_i : \mC_{X_i} \to \mC_i$, $\pi_i (\rho) := \rho_i$, such that
\begin{equation}
\label{eq_clutch}
\pi_i \circ \pi_j^{-1} = \alpha_{ij}
\ ;
\end{equation}
to be concise, in the previous equality we identified $\pi_i$ (resp. $\pi_j^{-1}$) with the restriction on $\mC_{X_{ij}}$ (resp. $\mC_{j,X_{ij}}$).

\section{Tensor \sC categories.}
\label{ssec_tens}

Let $\mT$ be a unital \sC category. As customary, we call {\em tensor product} an associative, unital \sC bifunctor $\otimes : \mT \times \mT \to \mT$ admitting an identity object $\iota \in \obt$, in such a way that $\iota \otimes \rho = \rho \otimes \iota = \rho$, $\rho \in \obt$. In such a case we say that $\mT$ is a {\em tensor} \sC category, and use the notation $(\mT , \otimes , \iota)$. For brevity, we adopt the notation $\rho \sigma := \rho \otimes \sigma$, so that if $t \in (\rs)$, $t' \in (\rsp)$, then $t \otimes t' \in (\rho \rho' , \sigma \sigma')$. Thus, if we pick $\rho' = \sigma' = \iota$, then we obtain that every space of arrows $(\rs)$ is a $(\ii)$-bimodule. In particular, $t \otimes 1_\iota = 1_\iota \otimes t = t$, $t \in (\rs)$. Note that $(\ii)$ is an Abelian \sC algebra with identity $1_\iota$. For basic properties of tensor \sC categories, we will refer to \cite{DR89,LR97}.
For every $\rho \in \obt$, we denote by $\wa \rho$ the tensor \sC category with objects $\rho^r$, $r \in \bN$, and arrows $(\rhors)$; we call $\rho$ the {\em generating object} of $\wa \rho$.

Duals of locally compact (quantum) groups are well-known examples of tensor \sC categories. Another example, which will be of particular interest in the present paper, is the category ${\bf vect}(X)$, with $X$ compact (Sec.\ref{ss_bun}). The tensor product on ${\bf vect}(X)$ is defined as in \cite[\S 1]{Ati},\cite[\S I.4]{Kar}, and will be denoted by $\otimes_X$ (no confusion should arise with the $C(X)$-algebra tensor product); note that $(\ii) = C(X)$.

\begin{ex}
\label{ex_end}
Let $\mA$ be a \sC algebra. Then, the set $\tend$ is endowed with a natural structure of tensor \sC category: the arrows are given by elements of the intertwiner spaces $(\rs) :=$ $\left\{ t \in \mA : t \rho(a) = \sigma (a) t , a \in \mA \right\}$, $\rs \in \tend$, and the tensor structure is defined by composition of endomorphisms (see \cite[\S 1]{DR89A}).
\end{ex}

\subsubsection{DR-dynamical systems.}
\label{ssec_dr}
Let $\rho \in \obt$. Then, a universal \sC dynamical system is associated with $\rho$, in the following way (see \cite[\S 4]{DR89} for details). For every $r,s \in \bN$, we consider the inductive limit Banach space $\oro^{r-s}$ associated with the sequence $i_{r,s} :$ $(\rhors) \hra$ $( \rho^{r+1} \rho^{s+1} )$, $i_{r,s} (t) :=$ $t \otimes 1_\rho$, $r \in \bN$; in particular, $\oro^0$ is the inductive limit \sC algebra associated with the sequence obtained for $r=s$. Composition of arrows and involution induce a *-algebra structure on $^0\oro := \sum_{k \in \bZ}^\oplus \oro^k$. Now, there exists a unique \sC norm on $^0\oro$ such that the {\em circle action}
\begin{equation}
\label{def_ca}
\wa z (t) := z^k t
\ , \ \
z \in \bT 
\ ,\ 
t \in \oro^k 
\ , \
k \in \bZ
\end{equation}
\noindent extends to an (isometric) automorphic action. We denote by $\oro$ the corresponding \sC algebra, called the {\em Doplicher-Roberts algebra} associated with $\rho$ ({\em DR-algebra}, in the sequel). A {\em canonical endomorphism} $\rho_* \in {\bf end} \oro$ can be defined,
\begin{equation}
\label{def_ce}
\rho_* (t) := 1_\rho \otimes t \ , \ \ t \in ( \rhors ) \ .
\end{equation}
\noindent We denote by $( \oro , \bT , \rho_* )$ the so-constructed \sC dynamical system. By construction, $\rho_* \circ \wa z = \wa z \circ \rho_*$, $z \in \bT$; moreover, $(\rhors) \subseteq ( \rho_*^r , \rho_*^s )$, $r , s \in \bN$, so that there is an inclusion of tensor \sC categories $\wa \rho \hra \wa \rho_*$. We say that $\rho$ is {\em amenable} if $\wa \rho = \wa \rho_*$, i.e. $(\rhors) = ( \rho_*^r , \rho_*^s )$, $r , s \in \bN$ (see \cite[\S 5]{LR97}). The above construction is universal, in the following sense: if $\pi : \wa \rho \to \wa \sigma$ is functor of tensor \sC categories, then there exists a \sC algebra morphism 
\begin{equation}
\label{def_um}
\wa \pi : \oro \to \mO_\sigma
\ , \ \
\wa \pi (t) := \pi (t) \ , \ t \in (\rhors) 
\end{equation}
\noindent such that $\sigma_* \circ \wa \pi = \wa \pi \circ \rho_*$ and $\wa z \circ \wa \pi = \wa \pi \circ \wa z$, $z \in \bT$. If $\pi$ is a \sC epifunctor, then $\wa \pi ( \rho_*^r , \rho_*^s ) \subseteq$ $( \sigma_*^r , \sigma_*^s )$; if $\sigma$ is amenable, then $\wa \pi ( \rho_*^r , \rho_*^s ) =$ $( \sigma_*^r , \sigma_*^s )$ and $\wa \pi$ is a \sC epimorphism.

\begin{ex}
\label{ex_od}
Let ${\bf hilb}$ denote the tensor \sC category of finite dimensional Hilbert spaces; if $\rho \in {\bf obj \ hilb}$, then $\oro$ is the {\em Cuntz algebra} $\mO_d$, where $d$ is the rank of $\rho$ (\cite{Cun77}). If $\mE \to X$ is a rank $d$ vector bundle, then the associated DR-algebra is the Cuntz-Pimsner algebra $\coe$ associated with the module of continuous sections of $\mE$. It turns out that $\coe$ is a locally trivial bundle with fibre $\mO_d$ (see \cite[\S 4]{Vas}, \cite{Vas05,Vas05B}).
\end{ex}

\subsubsection{Tensor categories as $C(X)$-categories.}
Let $( \mT , \otimes , \iota )$ be a tensor \sC category. We define $X^\iota$ as the spectrum of the Abelian \sC algebra
\begin{equation}
\label{def_cxi}
\left\{ 
f \in ( \ii ) : 
f \otimes 1_\rho = 1_\rho \otimes f \in (\rr)
\ , \
\rho \in \obt
\right\} \ ;
\end{equation}
\noindent note that (\ref{def_cxi}) has identity $1_\iota$, thus $X^\iota$ is compact. In the sequel, we will identify $C(X^\iota)$ with the \sC algebra (\ref{def_cxi}).

\begin{prop}
\label{prop_tens_cx}
Every tensor \sC category $( \mT , \otimes , \iota)$ has a natural structure of tensor $C (X^\iota)$-category, i.e. $\otimes$ is a $C(X^\iota)$-bifunctor.
\end{prop}

\begin{proof}
For every $\rho \in \obt$, we define the map $i_\rho : C (X^\iota) \to (\rr)$, $i_\rho (f) := f \otimes 1_\rho$. Since $i_\rho(1_\iota) = 1_\iota \otimes 1_\rho = 1_\rho$, it is clear that $i_\rho$ is a non-degenerate \sC algebra morphism. If $t \in (\rs)$, then $( 1_\sigma \otimes f ) \circ t = t \otimes f = t \circ ( 1_\rho \otimes f )$; thus $i_\sigma (f) \circ t = t \circ i_\rho (f)$, and $\mT$ is a $C (X^\iota)$-category. The fact that $\otimes$ preserves the $C(X^\iota)$-action follows from the obvious identities
\[
( 1_{\sigma \sigma'} \otimes f ) \circ (t \otimes t') =
(t \otimes f \otimes t') =
(t \otimes t' ) \otimes f =
( t \otimes t' ) \circ ( 1_{\rho \rho'} \otimes f ) \ ,
\]
\noindent $f \in C(X^\iota)$, $t \in (\rs)$, $t' \in (\rsp)$. 
\end{proof}

\begin{rem}
\label{rem_rsu}
According to the notation of Sec.\ref{sec_cx_cat}, in the sequel we will write $ft := f \otimes t$, $f \in C(X^\iota)$, $t \in (\rs)$. Let now $U \subseteq X^\iota$ be an open set; if $f \in C_U(X^\iota)$, $t \in (\rs)$, $t' \in (\rsp)$, then $ft \in (\rs)_U$, and it is clear that $ft \otimes t' \in ( \rho \rho' , \sigma \sigma' )_U$.
\end{rem}

It follows from Prop.\ref{prop_cx_cat} that for every $x \in X^\iota$ there exists a fibre functor $\pi_x : \mT \to \mT_x$. A structure of tensor \sC category is defined on $\mT_x$, by assigning
\begin{equation}
\label{def_tensx}
\left\{
\begin{array}{ll}
\rho_x \sigma_x := (\rho \sigma)_x \\
\pi_x (t) \otimes_x \pi_x(t') := \pi_x ( t \otimes t' )
\end{array}
\right.
\end{equation}
\noindent It is easy to prove that $\otimes_x$ is well-defined; in fact, if $v \in (\rs)$, $v' \in ( \rsp )$ and $\pi_x (v) = \pi_x(t)$, $\pi_x (v') = \pi_x (t')$ (i.e., $t-v \in (\rs)_{X - \left\{ x \right\}}$, $t'-v' \in (\rsp)_{X - \left\{ x \right\}}$), then $t-v = fw$, $t'-v' = f'w'$ for some $f,f' \in C_x(X)$, so that $t \otimes t' - v \otimes v'$ $=$ $fw \otimes v'$ $+$ $t \otimes f'w'$, and $\pi_x (t \otimes t') - \pi_x (v \otimes v') = 0$. This also proves that the fibre functors $\pi_x$, $x \in X$, preserve the tensor product. More generally, the above argument applies for the restriction $\mT_W$ over a closed $W \subset X^\iota$: a tensor product $\otimes_W : \mT_W \times \mT_W \to$ $\mT_W$ is defined, in such a way that the restriction epimorphism $\eta_W : \mT \to \mT_W$ preserves the tensor product.
We say that $\mT$ is a {\em continuous bundle of tensor \sC categories} if it is a continuous bundle w.r.t. the above $C(X^\iota)$-category structure.
As a consequence of the above considerations, we obtain a simple result on the structure of the DR-algebra associated with an object $\rho \in \obt$. The proof is trivial, therefore it will be omitted; note that we make the standard assumption that the map $t \mapsto t \otimes 1_\rho$, $t \in (\rhors)$, $r,s \in \bN$, is isometric.

\begin{prop}
\label{main_thm}
Let $\mT$ be a tensor \sC category (resp. a continuous bundle of tensor \sC categories), $\rho \in \obt$. Then, $\oro$ is a $C(X^\iota)$-algebra. In particular, if $\wa \rho$ is a (locally trivial) continuous bundle, then $\oro$ is a (locally trivial) continuous bundle of \sC algebras.
\end{prop}

A tensor \sC category $(\mT , \otimes , \iota)$ is said to be {\em symmetric} if for every $\rs \in \obt$ there exists a unitary $\eps (\rs) \in ( \rho \sigma , \sigma \rho )$ implementing the flip
\begin{equation}
\label{eq_flip}
\eps (\sigma , \sigma') \circ (t \otimes t') 
=
(t' \otimes t) \circ \eps (\rho , \rho') \ ,
\end{equation}
\noindent $t \in (\rs)$, $t' \in (\rsp)$, and satisfying certain natural relations (see \cite[\S 1]{DR89}). In essence, the notion of symmetry expresses commutativity of the tensor product. It is well-known that duals of locally compact groups are symmetric tensor \sC categories.

\begin{rem}
\label{rem_sym}
Let $( \mT , \otimes , \iota , \eps )$ be a symmetric tensor \sC category. By definition, $C(X^\iota) \subseteq (\ii)$. On the converse, if $z \in (\ii)$, $\rho \in \obt$, then by (\ref{eq_flip}) we find $z \otimes 1_\rho = \eps ( \rho , \iota ) \circ (1_\rho \otimes z) \circ \eps ( \iota , \rho )$; since $\eps ( \rho , \iota ) = \eps ( \iota , \rho ) = 1_\rho$, we conclude $z \in C(X^\iota)$ and $(\ii) =$ $C(X^\iota)$. For every $x \in X^\iota$, we denote by $\pi_x : \mT \to \mT_x$ the fibre epifunctor associated with $\mT$ as a $C(X^\iota)$-category. By (\ref{def_tensx}), every $\pi_x$ preserves the tensor product; moreover, by defining $\eps_x ( \rho_x , \sigma_x ) := \pi_x ( \eps (\rs) )$, we obtain that every $\mT_x$ is symmetric. By definition, the fibre functor $\pi_x$ preserves the symmetry. 
\end{rem}

Let $\rho \in \obt$. By \cite[Appendix]{DR89}, for every $r \in \bN$ there is a unitary representation of the permutation group of $r$ objects, which we denote by $\bP_r$:
\begin{equation}
\label{eq_def_epsrho}
\eps_\rho : \bP_r \to (\rho^r ,\rho^r) 
\ , \ \ 
p \mapsto \eps_\rho(p)
\ .
\end{equation}
\noindent For every $r \in \bN$, we introduce the antisymmetric projection
\[
P_{\rho,\eps,r}
:= 
\frac{1}{r!} \sum_{p \in \bP_r} {\mathrm{sign}} (p) \ \eps_\rho (p) \ . 
\]
\noindent Let $\rho_* \in {\bf end} \oro$ be the canonical endomorphism. It follows from (\ref{eq_flip}) that $\rho_*$ is "approximately inner", in the sense that
\begin{equation}
\label{eq_rhoeps}
\rho_* (t) = 1_\rho \otimes t = 
\eps_\rho (s,1) \ t \ \eps_\rho (1,r) 
\ , \ \ 
t \in (\rhors) \ 
\end{equation}
\noindent (recall that $t$ is identified with $t \otimes 1_\rho$ in $\oro$). For brevity (and coherence with (\ref{eq_def_epsrho})), in (\ref{eq_rhoeps}) we used the notation $\eps_\rho(r,s) :=$ $\eps (\rhors)$, $r,s \in \bN$.
Let $( \mT , \otimes , \iota , \eps )$ be a symmetric tensor \sC category. We denote by ${\bf aut}_\eps \mT$ the set of \sC autofunctors $\alpha$ of $\mT$ satisfying
\begin{equation}
\label{def_autr}
\alpha (\rho) = \rho \ , \ \
\alpha ( t \otimes t' ) = \alpha (t) \otimes \alpha (t')  \ , \ \
\alpha ( \eps (\rs) )  = \eps ( \rs ) \ ,
\end{equation}
\noindent $t \in (\rs)$, $t' \in (\rsp)$, $\eps (\rs) \in ( \rho \sigma , \sigma \rho )$. 
In particular, let us consider the case $\mT = \wa \rho$ for some object $\rho$. We note that in order to obtain the third of (\ref{def_autr}) it suffices to require $\alpha (\eps_\rho (1,1)) =$ $\eps_\rho (1,1)$. Moreover, $\wa z \in {\bf aut}_\eps \wa \rho$ for every $z \in \bT$, where $\wa z$ is defined by (\ref{def_ca}). We denote by ${\bf aut}_\eps \oro \subset {\bf aut} \oro$ the closed group of automorphisms commuting with $\rho_*$, and leaving every $\eps_\rho (r,s)$ fixed, $r,s \in \bN$. If $\rho$ is amenable in the sense of Sec.\ref{ssec_dr}, then by universality of $\oro$ we obtain an identification
\begin{equation}
\label{eq_aut}
{\bf aut}_\eps \wa \rho = {\bf aut}_\eps \oro \ .
\end{equation}
\noindent In fact, every \sC autofunctor $\alpha \in {\bf aut}_\eps \wa \rho$ defines an automorphism of $\oro$ as in (\ref{def_um}); on the converse, every $\beta \in {\bf aut}_\eps \oro$ defines by amenability a \sC autofunctor of $\wa \rho$. In the sequel we will regard ${\bf aut}_\eps \wa \rho$ as a topological group, endowed with the pointwise-convergence topology defined on ${\bf aut}_\eps \oro$. 
The following notion will play an important role in the sequel.
\begin{defn}
\label{def_ef}
Let $( \mT , \otimes  , \iota , \eps )$ be a symmetric tensor \sC category. An {\bf embedding functor} is a \sC monofunctor $i : \mT \hra {\bf vect} (X^\iota)$, preserving tensor product and symmetry.
\end{defn}

The following definition appeared in \cite[\S 2]{DR89}. A tensor \sC category $(\mT , \otimes , \iota)$ has {\em conjugates} if for every $\rho \in \obt$ there exists $\ovl \rho \in \obt$ with arrows $R \in ( \iota , \ovl \rho \rho )$, $\ovl R \in$$( \iota , \rho \ovl \rho )$, such that the {\em conjugate equations} hold:
\begin{equation}
\label{eq_conj}
({\ovl R}^* \otimes 1_\rho) \circ (1_\rho \otimes R) = 1_\rho
\ , \ \
(R^* \otimes 1_{\ovl \rho}) \circ (1_{\ovl \rho} \otimes \ovl R) = 1_{\ovl \rho}
\ .
\end{equation}
\noindent The notion of conjugation is deeply related with the one of dimension. Let us define
\begin{equation}
\label{def_dim}
d(\rho) = R^* \circ R \in (\ii) \ , \ d(\rho) > 0 \ ;
\end{equation}
\noindent then, it turns out that $d(\rho)$ is invariant w.r.t. unitary equivalence and conjugation (\cite[\S 2]{DR89}). If $\mT$ is symmetric, then (with some assumptions) $d(\rho) = d 1_\iota$ for some $d \in \bN$. In the case $\mT = {\bf hilb}$, then $d$ is the dimension in the sense of Hilbert spaces. For further investigations about the notion of dimension, we refer the reader to \cite{LR97,KPW04}.

\begin{ex}[Group duality in the Cuntz algebra]
\label{ss_cog}
The following construction appeared in \cite{DR87}. Let $\bH_d$ be the canonical rank $d$ Hilbert space. For every $r,s \in \bN$, we denote by $(\hrs)$ the Banach space of linear maps from $\bH^r$ into $\bH^s$; it is clear that $(\hrs)$ is isomorphic to the matrix space $\bM_{d^r , d^s}$. The category ${\wa \bH}_d$ with objects $\bH_d^r$, $r \in \bN$, and arrows $(\hrs)$ is a tensor \sC category, when endowed with the usual matrix tensor product. ${\wa \bH}_d$ is symmetric: we have the operators $\theta (r,s) \in ( \bH_d^{r+s} , \bH_d^{r+s} )$, $\theta (r,s) v \otimes v' := v' \otimes v$, $v \in \bH_d^r$, $v' \in \bH_d^s$. Moreover ${\bf aut}_\theta {\wa \bH}_d = {\bf aut}_\theta \mO_d \simeq \ud$ (\cite[Cor.3.3 , Lemma 3.6]{DR87}). The DR-dynamical system associated with $\bH_d$ is given by $( \mO_d , \bT , \sigma_d )$, where $\mO_d$ is the Cuntz algebra endowed with the {\em gauge action} $\bT \to {\bf aut} \mO_d$ and the {\em canonical endomorphism} $\sigma_d \in {\bf end} \mO_d$.
Let $G \subseteq \ud$ be a closed group. Then, every tensor power $\bH_d^r$, $r \in \bN$, is a $G$-module w.r.t. the natural action $g , \psi \mapsto g_r \psi$, $g \in G$, $g_r := g \otimes \ldots \otimes g$, $\psi \in \bH_d^r$. We denote by $\wa G$ the subcategory of ${\wa \bH}_d$ with arrows the spaces of equivariant maps
\begin{equation}
\label{eq_hrsg}
( \hrs )_G :=
\left\{ 
t \in ( \hrs )  :  g_s \circ t = t \circ g_r  \ , \ g \in G
\right\} \ .
\end{equation}
\noindent It is well-known that the dual of $G$ (i.e., the category of finite dimensional representations) is recovered by extending $\wa G$ w.r.t. direct sums. An automorphic action $G \to {\bf aut} \mO_d$ can be constructed, by defining
\begin{equation}
\label{def_cag}
\wa g \in {\bf aut}\mO_d 
\ : \ 
\wa g (t) := g_s \circ t \circ g_r^* \  ,  \ g \in G , t \in ( \hrs ) \ .
\end{equation}
\noindent Thus, every $( \hrs )_G$ coincides with the fixed-point space w.r.t. (\ref{def_cag}). We denote by $\mO_G$ the \sC subalgebra of $\mO_d$ generated by $(\hrs)_G$, $r,s \in \bN$; it turns out that $\mO_G$ is the fixed-point algebra w.r.t. (\ref{def_cag}). Since $\sigma_d \circ \wa g = \wa g \circ \sigma_d$, $g \in G$, we obtain that $\sigma_d$ restricts to an endomorphism $\sigma_G \in {\bf end} \mO_G$. It turns out that $( \hrs )_G = ( \sigma_G^r , \sigma_G^s )$, $r,s \in \bN$; in other terms, we find $\wa G = \wa \sigma_G$ (in particular, for $G = \left\{ 1 \right\}$, we have ${\wa \bH}_d = {\wa \sigma}_d$). Moreover, the stabilizer ${\bf aut}_{\mO_G} \mO_d$ of $\mO_G$ in $\mO_d$ is isomorphic to $G$ {\em via} (\ref{def_cag}). Thus, $G$ and the associated tensor \sC category $\wa G$ are recovered by properties of the \sC dynamical system $( \mO_G , \sigma_G )$, together with the inclusion $\mO_G \hra \mO_d$.
\end{ex}

Let $(\mT , \otimes , \iota , \eps)$ be a symmetric tensor \sC category. $\rho \in \obt$ is said to be a {\em special object} if there exists $d \in \bN$, $V \in (\iota , \rho^d)$ such that 
\begin{equation}
\label{eq_det}
\left\{
\begin{array}{ll}
V^* \circ V = 1_\iota \ , \ \ V \circ V^* = P_{\rho,\eps,d} \\
(V^* \otimes 1_\rho ) \circ ( 1_\rho \otimes V ) = (-1)^{d-1} d^{-1} 1_\rho
\ .
\end{array}
\right.
\end{equation}
\noindent It turns out that $V$ is unique up to multiplication by elements of $\bT$ (see \cite[\S 3]{DR89}). In the sequel, we will denote the above data by $(\rho , d , V)$. It is proved that the integer $d$ coincides with the dimension of $\rho$ (\cite[Lemma 3.6]{DR89}). 
Special objects play a pivotal role in the Doplicher-Roberts duality theory. In fact, tensor categories generated by special objects are in one-to-one correspondence with duals of compact Lie groups, in the case $(\ii) \simeq \bC$. Since in the sequel we will make use of such concepts at a detailed level, in the next Lemma an abstract duality is summarized for special objects.

\begin{lem}
\label{thm_spec}
Let $( \rho , d , V )$ be a special object, and suppose $(\ii) \simeq \bC$. Then, the following properties are satisfied:
\begin{enumerate}
\item  $\rho$ is amenable, and there is a \sC monomorphism $i : ( \oro , \rho_* ) \hra ( \mO_d , \sigma_d )$;
\item  the above monomorphism defines an embedding functor $i_* : \wa \rho \hra \wa \bH_d$, in the sense of Def.\ref{def_ef};
\item  there is a closed group $G \subseteq \sud$ such that $i_* (\wa \rho) = \wa G$, so that $\wa G$, $\wa \rho$ are isomorphic as symmetric tensor \sC categories;
\item  $i$ restricts to an isomorphism $( \oro , \rho_* ) \simeq ( \mO_G , \sigma_G )$;
\item  $G$ is isomorphic to the stabilizer of $\mO_G$ in $\mO_d$ w.r.t. the action (\ref{def_cag});

\item  if $\alpha \in {\bf aut}_\eps \wa \rho$, then there is a unitary $u \in \ud$ such that $i \circ \alpha = \wa u \circ i$ (where $\wa u \in {\bf aut}_\theta \mO_d$ is defined as in (\ref{def_cag})).
\end{enumerate}
\noindent 
\end{lem}

\begin{proof}
Point 1 is \cite[Lemma 4.14]{DR89}. Points 2,3,4,5 are proved in \cite[Thm.4.17]{DR89}. About point 6, let $\alpha \in {\bf aut}_\eps \wa \rho$. Since $\rho_* (V) = \eps_\rho ( d,1 ) V$ (recall (\ref{eq_rhoeps})), we obtain 
\[
\rho_* ( V^* \alpha ( V ) ) =
V^* \eps_\rho ( 1,d  ) \alpha ( \eps_\rho ( d,1 ) V ) =
V^* \alpha (V) \ . 
\]
\noindent Thus, $V^* \alpha (V)$ is $\rho_*$-invariant, and $V^* \alpha(V) \in (\ii) = \bC$ (see \cite[Lemma 4.13]{DR89}). Now, $\alpha ( P_{\rho,\eps,d}  ) = P_{\rho,\eps,d}$, so that 
\[
1 = \left\| V \right\|^2 
\geq 
\left\| V^* \alpha (V) \right\| 
=
\left\| V \right\| \left\| V^* \alpha (V) \right\| \left\| \alpha (V^*) \right\|
\geq
\left\| P_{\rho,\eps,d}^2 \right\|
= 1
\ , 
\]
\noindent we conclude that $V^* \alpha (V)  = z_\rho$, $z_\rho \in \bT$. Thus, 
\[
\alpha (V) = P_{\rho,\eps,d} \alpha (V) = VV^* \alpha(V) = z_\rho V \ .
\]
\noindent Let now $\lambda_\rho \in \bT$, $\lambda_\rho^d = \ovl z_\rho$. We define the automorphism $\beta := \wa \lambda_\rho \circ \alpha \in {\bf aut}_\eps \wa \rho$, in such a way that $\beta (V) = \wa \lambda_\rho \circ \alpha (V) = \lambda_\rho^d z_\rho V = V$. By \cite[Cor.4.9(d)]{DR89A}, we obtain that there is $v \in \sud$ such that $i \circ \beta = \wa v \circ i$. By defining $u := \ovl \lambda_\rho v \in \ud$, we obtain $\wa u \circ i = i \circ \alpha$. 
\end{proof}

\subsubsection{Twisted special objects.}
Lemma \ref{thm_spec} implies that tensor \sC categories generated by special objects correspond to duals of compact Lie groups $G \subseteq \sud$, in the case $(\ii) \simeq \bC$. In particular, by performing the embedding $\wa \rho \hra \wa \bH_d$, we obtain that $V$ corresponds to a normalized generator of the totally antisymmetric tensor product $\wedge^d \bH_d$. 
In the case in which $(\ii)$ is nontrivial, the notion of special object is too narrow to cover all the interesting cases. For example, let us consider the category ${\bf vect} (X)$, where $X$ is a compact Hausdorff space: if $\mE \to X$ is a rank $d$ vector bundle, then it is well-known that the totally antisymmetric tensor product $\wedge^d \mE$ is a line bundle, in general non-trivial. The cohomological obstruction to get triviality of $\wedge^d \mE$ is encoded by the first Chern class of $\mE$. In the particular case in which $\wedge^d \mE$ is trivial, then we may find a (normalized) section, which plays the right role in the definition of special object. Since our duality theory will be modeled on ${\bf vect}(X)$, it becomes natural to generalize the notion of special object, in such a way to include the above-described case.

Let $\mT$ be a symmetric tensor \sC category. A {\em twisted special object} is given by a triple $(  \rho , d , \mV  )$, where $\rho \in \obt$, $d \in \bN$, and $\mV \subseteq ( \iota , \rho^d )$ is a closed subspace such that for every $V, V' \in \mV$, $f \in (\ii)$,
\begin{equation}
\label{eq_tso}
\left\{
\begin{array}{ll}
f \otimes V = V \otimes f \in \mV    \\
(V^* \otimes 1_\rho ) \circ ( 1_\rho \otimes V' ) \ = \ 
(-1)^{d-1} d^{-1} ( V^* \circ V' ) \otimes 1_\rho   \\
{\mV}{\mV}^*  \ := \ 
{\mathrm{span}}\left\{ V' \circ  V^*  :  V , V' \in \mV_\rho \right\} \ = \ 
(\ii) \otimes P_{\rho,\eps,d} \ .
\end{array}
\right.
\end{equation}
\noindent Note that $V^* V' \in (\ii)$, thus $\mV$ is endowed with a natural structure of right Hilbert $(\ii)$-module. It is clear that ${\mV}{\mV}^*$ may be identified with the \sC algebra of compact, right $(\ii)$-module operators on $\mV$, and that $P_{\rho,\eps,d}$ is the identity of $\mV$. The map $i : (\ii) \to K(\mV)$, $i(f) := f \otimes P_{\rho,\eps,d}$ defines a left $(\ii)$-action on $\mV$; note that (\ref{eq_tso}(3)) implies that $i$ is surjective. From (\ref{eq_tso}(3)) we also obtain that $\mV$ is full, in fact $\left\| V^* \circ V \right\| = $ $\left\| V \circ V^* \right\| = $ $\left\| f_V \otimes P_{\rho,\eps,d} \right\|$, $V \in \mV$, $f_V \in (\ii)$, so that every positive element of $(\ii)$ appears as the square of the norm of an element of $\mV$. Moreover, $i(f) V =$ $f \otimes V = 0$, $V \in \mV$, implies $f \otimes V^* V = 0$, i.e. $ff_V = 0$, $f_V :=$ $V^* V \in$ $(\ii)$; since at varying of $V$ in $\mV$ we obtain all the positive elements of $(\ii)$, we find that $f = 0$, and the left action $i$ is injective. We conclude that $\mV$ is a bi-Hilbertian $(\ii)$-bimodule.

\begin{lem}
If existing, the Hilbert $(\ii)$-bimodule $\mV$ is unique. Moreover, there exists a unique up-to-isomorphism line bundle $\mL_\rho \to X^\iota$ such that $\mV$ is isomorphic to the module of continuous sections of $\mL_\rho$.
\end{lem}

\begin{proof}
The generalised Serre-Swan theorem proved in \cite{DG} implies that $\mV$ is the module of continuous sections of a Hilbert bundle $\mL_\rho \to X^\iota$; since $K(\mV) \simeq C(X^\iota)$, we conclude that $\mL_\rho$ has rank one, i.e. it is a line bundle. 
We now pass to prove the unicity; we consider a finite set of generators $\left\{ V_i \right\} \subset \mV$, so that $\sum_i V_i \circ V_i^* = P_{\rho,\eps,d}$, $\sum_i V_i^* \circ V_i = 1_\iota$. If $\mW \subseteq ( \iota , \rho^d )$ is another $(\ii)$-bimodule satisfying (\ref{eq_tso}) and $W \in \mW$, then $W = P_{\rho,\eps,d} \circ W$. This implies $W = \sum_i V_i \circ (V_i^* \circ W)$, with $V_i^* \circ W \in (\ii)$, thus $W \in \mV$. This proves $\mW = \mV$. 
\end{proof}

The set of transition maps associated with $\mL_\rho \to X^\iota$ defines a $\bT$-cocycle in $H^1(X^\iota,\bT)$. The isomorphism $H^1(X^\iota,\bT) \simeq H^2(X^\iota,\bZ)$ allows one to give the following definition.

\begin{defn}
\label{def_cc_tso}
The {\bf Chern class} of a twisted special object $( \rho , d , \mV )$ is the unique element $c (\rho) \in H^2(X^\iota,\bZ)$ associated with $\mL_\rho$.
\end{defn}

By definition, usual special objects are exactly those with Chern class $c(\rho) = 0$; in such a case, $\mV$ is generated as a $(\ii)$-module by an isometry $V \in (\iota , \rho^d)$. Let now $W \subseteq X^\iota$; we denote by $\eta_W : \mT \to \mT_W$ the restriction epimorphism; by using the argument of Rem.\ref{rem_sym}, for every $\rs \in \obt$ we define $\eps_W (\rs) := \eta_W (\eps(\rs))$, and conclude that $(\mT_W ,\otimes_W , \eta_W(\iota) , \eps_W )$ is symmetric.

\begin{cor}
For every $x \in X^\iota$, there exists a closed neighborhood $W \subseteq X^\iota$, $W \ni x$, such that $\eta_W (\rho)$ is a special object.
\end{cor}

\begin{proof}
It suffices to pick a closed neighborhood $W$ trivializing $\mL_\rho$: this implies the existence of a normalized section of $\mL_\rho |_W$, which appears as an element $V_W \in \eta_W (\mV) \subseteq$ $\eta_W (\iota , \rho^d)$ such that $V_W^* \circ V_W = \eta_W (1_\iota)$, $V_W \circ V_W^* = \eta_W(P_{\rho,\eps,d})$. 
\end{proof}

Let $\mA$ be a \sC algebra. Twisted special objects in ${\bf end} \mA$ have been studied in \cite[\S 4]{Vas05A}, and have been called {\em special endomorphisms}.

\begin{lem}
\label{lem_oro_tso}
Let $( \rho , d , \mV )$ be a twisted special object. Then $\rho$ is amenable, and $\wa \rho$ is a continuous bundle of tensor \sC categories with base space $X^\iota$. For every $x \in X^\iota$, there exists a compact Lie group $G(x) \subseteq \sud$ such that the fibre $\wa \rho_x$ is isomorphic to $\wa{G(x)}$.
\end{lem}

\begin{proof}
As first, we prove that $\rho$ is amenable. Let $\left\{ X_i \right\}$ be an open cover trivializing $\mL_\rho$, with a subordinate partition of unity $\left\{ \lambda_i \right\}$. We denote by $\wa \rho_i$ the restriction of $\wa \rho$ on the closure $\ovl X_i$, and by $\eta_i : \wa \rho \to \wa \rho_i$ the restriction epifunctor. Since $\mL_\rho |_{X_i}$ is trivial, we find that $( \rho_i , d , V_i )$ is a special object, where $V_i \in \eta_i (\mV) \subseteq$ $( \iota , \rho_i^d )$ is a suitable partial isometry. By \cite[Lemma 4.14]{DR89} we find that $\rho_i$ is amenable, i.e. $( \rho_{i,*}^r , \rho_{i,*}^s ) =$ $(\rho_i^r , \rho_i^s)$, $r,s \in \bN$. This implies that $\eta_i$ induces a \sC epimorphism $\wa \eta_i : ( \oro , \rho_* ) \to$ $( \mO_{\rho_i} , \rho_{i,*} )$. 
Let $t \in ( \rho_*^r , \rho_*^s )$; for every index $i$, we consider $\wa \eta_i (t) \in ( \rho_{i,*}^r , \rho_{i,*}^s )$. Since $\rho_i$ is amenable, we conclude that $\wa \eta_i (t) \in ( \rho_i^r , \rho_i^s )$. By construction, $\eta_i : (\rhors) \to$ $( \rho_i^r , \rho_i^s )$ is the restriction on $\ovl X_i$ of the continuous field of Banach spaces $(\rhors)$; by the Tietze theorem \cite[10.1.12]{Dix}, there is $t_i \in (\rhors)$ such that $\eta_i (t_i) = \wa \eta_i (t)$. Now, we have $t = \sum_i \lambda_i t_i$; since $\sum_i \lambda_i t_i \in (\rhors)$, we conclude that $t \in (\rhors)$, and $\rho$ is amenable.
We now prove that $\wa \rho$ is a continuous bundle. Since $\rho$ is a twisted special object, we find that $\rho_* \in$ ${\bf end} \oro$ is a quasi-special endomorphism in the sense of \cite[Def.4.10]{Vas05A}, with the additional property that the permutation symmetry \cite[Def.1.1,\S 4.1]{Vas05A} is satisfied. Thus, by applying \cite[Thm.5.1,Cor.5.2]{Vas05A}, we find that $\oro$ is a continuous bundle, with fibres isomorphic to $\mO_{G(x)}$, $x \in X$, where each $G(x) \subseteq \sud$ is a compact Lie group. The same results also imply that $\wa \rho$ is a continuous bundle with fibres $\wa{G(x)}$, $x \in X$. 
\end{proof}

Let $( \rho , d , \mV )$ be a twisted special object. Then, $\mV$ appears as a Hilbert $(\ii)$-bimodule in $\oro$, so that an inner endomorphism $\nu \in {\bf end} \oro$ is defined, in the sense of \cite[\S 3]{Vas05}. In explicit terms, if $\left\{ \psi_l \right\} \subset$ $\mV$ is a finite set of generators for $\mV$, then
\begin{equation}
\label{def_nu}
\nu (t) := \sum_l \psi_l t \psi_l^* 
\ , \ \
t \in \oro \ .
\end{equation}
\noindent Since $\wa z (\psi_l) =$ $z^d \psi_l$, $z \in \bT$, we find that $\nu$ commutes with the circle action, i.e. $\nu \circ \wa z =$ $\wa z \circ \nu$, $z \in \bT$. Thus, in particular we find $\nu (\oro^k) \subseteq$ $\oro^k$, $k \in \bZ$. 

\subsubsection{Continuity of tensor \sC categories.}
\label{ssec_sym_b}
The following notion is a generalization of the analogue in \cite[\S 3]{DR89}. We say that a symmetric tensor \sC category $( \mT , \otimes , \iota , \eps )$ is {\em T-specially directed} if every object $\rho' \in \obt$ is {\em dominated} by a twisted special object $\rho \in \obt$, i.e. there exist orthogonal partial isometries $S_i \in (\rho' , \rho^{n_i})$, $i = 1,2, \ldots , m$, with $1_{\rho'} = \sum_i S_i^* \circ S_i$. 

If $\mT$ has direct sums, subobjects and conjugates, then $\mT$ is T-specially directed, with the additional property that every $\rho'$ is dominated by a (non-twisted) special object (see the proof of \cite[Thm.3.4]{DR89}).

\begin{thm}
\label{thm_tcf}
Let $( \mT , \otimes , \iota , \eps )$ be a symmetric tensor \sC category. Then, $C(X^\iota)$ coincides with $( \ii )$. If $\mT$ is T-specially directed, then $\mT$ is a continuous bundle of symmetric tensor \sC categories over $X^\iota$. 
\end{thm}
\begin{proof}
In Rem.\ref{rem_sym}, we verified that $C(X^\iota) = (\ii)$, thus it remains to prove that $\mT$ is a continuous bundle. At this purpose, we consider $\rs \in \obt$, $t \in (\rs)$, and prove the continuity of the norm function $n_t (x) := \left\|  \pi_x (t)  \right\|$, $x \in X^\iota$. Since $n_t (x) =$ $\left\| \pi_x (t^* \circ t) \right\|^{1/2}$, with $t^* \circ t \in$ $( \rr )$, it suffices to verify the continuity of $n_t$ only for arrows that belong to $( \rr )$, $\rho \in \obt$.
By Lemma \ref{lem_oro_tso}, the norm function $n_t$ is continuous for every $t \in (\rhors)$, $r,s \in \bN$, where $\rho$ is a twisted special object. Since there is an obvious inclusion $\wa{\rho^n} \hra$ $\wa \rho$, the norm function remains continuous for arrows between tensor powers of twisted special objects.
This implies that if we consider the full \sC subcategory $\mT^{sp}$ of $\mT$ with objects tensor powers of twisted special objects, then $\mT^{sp}$ is a continuous bundle of \sC categories over $X^\iota$. By Lemma \ref{lem_ccx}, we conclude that $(\mT^{sp})_{s,+}$ is a continuous bundle of \sC categories. 
Since $\mT$ is $T$-specially directed, we find that every $\rho' \in \obt$ is an object of $(\mT^{sp})_{s,+}$, in fact $\rho'$ is the direct sum of subobjects of $\rho^{n_i}$, $i = 1 , \ldots , m$, where $\rho$ is a twisted special object.
We conclude that the norm function of $\rho'$ is continuous, and the theorem is proved. 
\end{proof}

\begin{rem}
\label{rem_2c}
An analogue of the previous theorem has been proved by P. Zito (\cite[\S 3]{Zit05}), in the setting of $2$-\sC categories with conjugates. In the above-cited result, no symmetry property of the tensor product is assumed. It is not difficult to prove that the bundle structure constructed by Zito coincides with the one of Thm.\ref{thm_tcf}, in the common case of symmetric tensor \sC categories with conjugates.
\end{rem}

\begin{rem}
\label{rem_fib_dual}
The fibre epifunctors $\pi_x$, $x \in X^\iota$, considered in the previous theorem preserve symmetry and tensor product, thus every fibre $\mT_x$ has conjugates; moreover, $\pi_x ( \ii ) \simeq \bC$ for every $x \in X$. By closing each $\mT_x$ w.r.t. subobjects and direct sums, we obtain a tensor \sC category satisfying the hypothesis of \cite[Thm.6.1]{DR89}, which turns out to be isomorphic to the dual of a compact group ${G(x)}$. In particular, for each pair $\rs \in \obt$, the fibres of the continuous field $(\rs)$ are the finite-dimensional vector spaces $( \bH_{\rho,x} , \bH_{\sigma,x} )_{{G(x)}}$ of ${G(x)}$-equivariant operators between Hilbert spaces $\bH_{\rho,x}$, $\bH_{\sigma,x}$, $x \in X^\iota$.
\end{rem}

\subsubsection{Local Triviality.}
In the next definition, we give a notion of local triviality for a symmetric tensor \sC category, compatible with the bundle structure of Thm.\ref{thm_tcf}. Let $( \mT_\bullet , \otimes_\bullet , \iota_\bullet , \eps_\bullet )$ be a symmetric tensor \sC category, $X$ a compact Hausdorff space; then, the constant bundle $X \mT_\bullet$ (Sec.\ref{sec_cx_cat}) has a natural structure of symmetric tensor \sC category,
\begin{equation}
\label{def_ts}
( t \otimes_\bullet^X t' ) (x) := t(x) \otimes_\bullet t'(x) 
\ , \ \ 
\eps_\bullet^X (\rs) := \eps_\bullet (\rs)  \ ,
\end{equation}
\noindent  $t,t' \in (\rs)^X$, $\rs \in {\bf obj} \ \mT_\bullet$. 
Let now $( \mT , \otimes , \iota , \eps )$ be a symmetric tensor \sC category. The fibres $\mT_x$, $x \in X^\iota$, are symmetric tensor \sC categories having the same set of objects as $\mT$, in such a way that the fibre epifunctors $\pi_x : \mT \to \mT_x$ induce the identity map $\obt \to {\bf obj} \ \mT_x \equiv \obt$ (Rem.\ref{rem_str_fibre}); moreover, $\pi_x$ preserves tensor product and symmetry (Rem.\ref{rem_sym}), and every $\mT_x$, $x \in X^\iota$, has simple identity object.

\begin{defn}
\label{def_slc}
Let $( \mT , \otimes , \iota , \eps )$ be a symmetric tensor \sC category, endowed with the natural $C(X^\iota)$-category structure. $\mT$ is said to be {\bf locally trivial} if it is locally trivial in the sense of Sec.\ref{ss_lt_bo}, with the additional property that the local charts preserve tensor product and symmetry.
\end{defn}

In explicit terms, there is a symmetric tensor \sC category $( \mT_\bullet , \otimes_\bullet , \iota_\bullet , \eps_\bullet )$ with simple identity object, and a cover $\left\{ X_i \subseteq X^\iota \right\}$ of closed neighborhoods with \sC isofunctors $\pi_i : \mT_{X_i} \to X_i \mT_\bullet$ satisfying
\begin{equation}
\label{def_lc}
\left\{
\begin{array}{ll}
\pi_i (\rho) = \pi_j (\rho)  \\
\pi_i ( t \otimes t' ) = \pi_i (t) \otimes_\bullet^{X_i} \pi_i (t')  \\
\pi_i ( \eps (\rs) )  = \eps_\bullet^{X_i} (  \rs  )
\end{array}
\right.
\end{equation}
\noindent where $\rs \in \obt$, $t \in (\rs)$, $t' \in (\rsp)$. If $X$ is a compact Hausdorff space, we denote by
\begin{equation}
\label{not_tens}
\sym ( X , \mT_\bullet )
\end{equation}
\noindent the set of isomorphism classes of locally trivial, symmetric tensor \sC categories with fibres isomorphic to $\mT_\bullet$, and such that the space of intertwiners of the identity object is isomorphic to $C(X)$.

Let $( \wa \rho_\bullet , \otimes_\bullet , \iota_\bullet , \eps_\bullet )$ be a symmetric tensor \sC category with generating object $\rho_\bullet$, such that $( \iota_\bullet , \iota_\bullet ) \simeq \bC$. If $( \mT , \otimes , \iota , \eps ) \in \sym ( X , \wa \rho_\bullet )$, then $\mT$ has the same objects as $\wa \rho_\bullet$; in order to avoid confusion, we will denote by $\rho$ the object of $\mT$ corresponding to $\rho_\bullet$, so that $\mT$ is generated by the tensor powers of $\rho$, i.e. $\mT = \wa \rho$.

In the following lines, we regard the topological group $\autrb$ (\ref{def_autr}) as the "structure group" for a locally trivial, symmetric tensor \sC category.
Let $K$ be a topological group. A $K$-cocycle is given by a pair $( \left\{ X_i \right\} , \left\{ g_{ij} \right\} )$, where $\left\{ X_i \right\}$ is a finite open cover of $X$ and $g_{ij} : X_{ij} \to K$ are continuous maps such that $g_{ij} (x) g_{jk}(x) =$ $g_{ik} (x)$, $x \in X_{ijk}$. We say that cocycles $( \left\{ X_i \right\} , \left\{ g_{ij} \right\} )$, $( \left\{ X'_l \right\} , \left\{ g'_{lm} \right\} )$ are equivalent if there are continuous maps $u_{il} : X_i \cap X'_l \to K$ such that $g_{ij} (x) u_{jm} (x) =$ $u_{il} (x) g'_{lm} (x)$, $x \in X_{ij} \cap X'_{lm}$. The set of equivalence classes of $K$-cocycles is denoted by $H^1(X,K)$. It is well-known that $H^1(X,K)$ classifies the principal $K$-bundles over $X$ (\cite[Chp.4]{Hus}).

\begin{lem}
\label{lem_lc}
Let $\rho_\bullet$ be amenable. Then, there is a one-to-one correspondence $\sym ( X , \wa \rho_\bullet ) \lra H^1 ( X , \autrb )$.
\end{lem}

\begin{proof}
Let $( \wa \rho , \otimes , \iota , \eps ) \in \sym ( X , \wa \rho_\bullet )$. Then, there are local charts $\pi_i : \wa \rho_{X_i} \to X_i \wa \rho_\bullet$, where $\left\{ X_i \subseteq X \right\}$ is a cover of closed neighborhoods. Let $X_{ij} \neq \emptyset$; then, every $\pi_i$ restricts in a natural way to a local chart $\pi_{i,ij} : \wa \rho_{X_{ij}} \to X_{ij} \wa \rho_\bullet$. We define $\alpha_{ij} : X_{ij} \wa \rho_\bullet \to X_{ij} \wa \rho_\bullet$, $\alpha_{ij} := \pi_{i,ij} \circ \pi_{j,ij}^{-1}$; by (\ref{def_lc}), we obtain
\[
\alpha_{ij} (t \otimes_\bullet^{X_{ij}} t') = 
\alpha_{ij} (t) \otimes_\bullet^{X_{ij}} \alpha_{ij} (t')  
\ , \ \
\alpha_{ij} ( \eps_\bullet^{X_{ij}} (r,s) ) = \eps_\bullet^{X_{ij}} ( r,s )
\ ,
\]
\noindent $t \in (\rho_\bullet^r , \rho_\bullet^s)_{X_{ij}}$, $t' \in (\rho_\bullet^{r'} , \rho_\bullet^{s'} )_{X_{ij}}$. Thus, $\alpha_{ij} \in {\bf aut}_{\eps_\bullet} (  X_{ij} \wa \rho_\bullet  )$. Now, the DR-algebra associated with $\rho_\bullet$ (regarded as an object of $X_{ij} \wa \rho_\bullet$) is the trivial field $C(X_{ij}) \otimes \mO_{\rho_\bullet}$; we denote by $\pi_x : C(X_{ij}) \otimes \mO_{\rho_\bullet} \to \mO_{\rho_\bullet}$ the evaluation epimorphism assigned for $x \in X_{ij}$. By (\ref{eq_aut}), we may regard $\alpha_{ij}$ as an element of ${\bf aut}_{\rho_\bullet , \eps_\bullet} ( C(X_{ij}) \otimes \mO_{\rho_\bullet} )$. In particular, $\alpha_{ij}$ is a $C(X_{ij})$-automorphism, thus we may regard $\alpha_{ij}$ as a continuous map
\[
\alpha_{ij} : X_{ij} \to {\bf aut}_{\rho_\bullet , \eps_\bullet} \mO_{\rho_\bullet}
\ , \ 
[\alpha_{ij} (x)] (t) := \pi_x \circ \alpha_{ij} (1_{ij} \otimes t)
\ ,
\]
\noindent where $t \in \mO_{\rho_\bullet}$, $x \in X_{ij}$, and $1_{ij} \in C(X_{ij})$ is the identity. Now, the obvious identity $\pi_{i,ijk} \circ \pi_{k,ijk}^{-1}  =$ $\pi_{i,ijk} \circ \pi_{j,ijk}^{-1} \circ \pi_{j,ijk} \circ \pi_{k,ijk}^{-1}$ implies $\alpha_{ik} = \alpha_{ij} \circ \alpha_{jk}$ (over $X_{ijk}$); thus, by applying (\ref{eq_aut}), we conclude that the set $\left\{ \alpha_{ij} \right\}$ defines an $\autrb$-cocycle. By choosing another set of local charts $\pi'_h : \wa \rho \to Y_h \wa \rho_\bullet$, we obtain a cocycle $\beta_{hk} := \pi'_{h,hk} \circ (\pi'_{k,hk})^{-1}$ which is equivalent to $\left\{ \alpha_{ij} \right\}$, in fact $\beta_{hk} = V_{hi} \circ \alpha_{ij} \circ V_{kj}^{-1}$, $V_{hi} : X_i \cap Y_h \to \autrb$, $V_{hi} := \pi'_h \circ \pi_i^{-1}$. If $( \wa {\rho'} , \otimes' , \iota' , \eps' )$ is a symmetric tensor \sC category with a \sC isofunctor $F : \wa \rho \to \wa {\rho'}$, then $\left\{ \pi_i \circ F \right\}$ is a set of local charts associated with $\wa{\rho'}$; thus, $\alpha_{ij} = ( \pi_i \circ F ) \circ ( F^{-1} \circ \pi_j^{-1} )$ is a cocycle associated with $\wa{\rho'}$. In other terms, we defined an injective map
\[
i_\bullet : \ \sym ( X , \wa \rho_\bullet  ) \hra H^1 (  X , \autrb ) \ .
\]
\noindent On the converse, let $\left\{ \alpha_{ij} \right\}$ be an $\autrb$-cocycle associated with a finite cover $\left\{ X_i \right\}$ of closed neighborhoods. As a preliminary remark, we consider the obvious structure of symmetric tensor \sC category (\ref{def_ts}) on the trivial bundle $W \wa \rho_\bullet$, $W \subseteq X$ closed. 
With such a structure, every $\alpha_{ij}$ defines a \sC autofunctor $\alpha_{ij} : X_{ij} \wa \rho_\bullet \to X_{ij} \wa \rho_\bullet$ preserving tensor product and symmetry, and such that $\alpha_{ij} (\rho_\bullet) = \rho_\bullet$. We consider the set of symmetric tensor \sC categories $\left\{ ( X_i \wa \rho_\bullet , \otimes_i , \iota_i , \eps_i ) \right\}$ (where $\otimes_i$, $\iota_i$, $\eps_i$ are defined according to (\ref{def_ts})), and the \sC category $\mT$ obtained by glueing every $X_i \wa \rho_\bullet$ w.r.t. the \sC isofunctors $\alpha_{ij}$'s. By construction, the objects of $\mT$ are the tensor powers of $\rho_\bullet$; the spaces of arrows of $\mT$ are defined according to (\ref{def_glue}). We define a tensor product on $\mT$, by posing $\rho_\bullet^r \rho_\bullet^s := \rho_\bullet^{r+s}$, $r,s \in \bN$, and
\begin{equation}
\label{def_tp}
( t_i )_i \otimes ( t'_i )_i
:=
( t_i \otimes_i t'_i )_i
\end{equation}
\noindent where the families $\left\{ t_i \in C(X_i , ( \rho_\bullet^{r} , \rho_\bullet^{s} )) \right\}$, $\left\{ t'_i \in C(X_i , ( \rho_\bullet^{r'} , \rho_\bullet^{s'} )) \right\}$ satisfy (\ref{def_glue}). Since each $\alpha_{ij}$ preserves the tensor product, we find 
\[
\begin{array}{ll}
\eta_{i,X_{ij}} (t_i \otimes t'_i) & =
\eta_{j,X_{ij}} (t_i) \otimes_i \eta_{j,X_{ij}} (t'_i)     \\ & =
\alpha_{ij} (\eta_{j,X_{ij}} (t_j)) 
\otimes_j 
\alpha_{ij} (\eta_{j,X_{ij}} (t'_j))     \\ & =
\alpha_{ij} \circ \eta_{j,X_{ij}} (t_j \otimes_j t'_j)
\ ,
\end{array}
\]
\noindent so that the l.h.s. of (\ref{def_tp}) is actually an arrow in $\mT$; thus, the tensor product on $\mT$ is well-defined. For the same reason, the operators
\[
\eps (r,s) := 
( \
1_{X_i} \otimes \eps_\bullet (\rho_\bullet^{r} , \rho_\bullet^{s}) 
\ )_i
\]
\noindent define a symmetry on $\mT$. We denote by $\mT := ( \wa \rho , \otimes , \iota , \eps )$ the symmetric tensor \sC category obtained in this way. By construction, $\wa \rho$ is equipped with a set of local charts $\pi_i : \wa \rho_{X_i} \to X_i \wa \rho_\bullet$ such that $\alpha_{ij} = \pi_{i,ij} \circ \pi_{j,ij}^{-1}$ (see (\ref{eq_clutch})), thus the map $i_\bullet$ is also surjective. 
\end{proof}

\begin{ex}[The Permutation category.]
\label{ex_pc}
Let $( \wa \rho_\bullet , \otimes_\bullet , \iota_\bullet , \eps_\bullet )$ be a symmetric tensor \sC category with $( \iota_\bullet , \iota_\bullet ) \simeq \bC$. We assume that $(\rho_\bullet^r , \rho_\bullet^s) = \left\{ 0 \right\}$ if $r \neq s \in \bN$, and that every $(\rho_\bullet^r , \rho_\bullet^r)$ is generated as a Banach space by the permutation operators ${\eps_\bullet}_{\rho_\bullet} (p)$, $p \in \bP_r$. It is clear that in this case the automorphism group $\autrb$ reduces to the identity. Let $d \in \bN$ denote the dimension of $\rho_\bullet$; by \cite[Lemma 2.17]{DR89} we find that $\wa \rho_\bullet$ is isomorphic as a symmetric tensor \sC category to $\wa \ud$ (see Ex.\ref{ss_cog}). Thus, $( \rho_\bullet^r , \rho_\bullet^r ) \simeq (\bH_d^r , \bH_d^r)_\ud$, $r \in \bN$. Let $\wa \rho \in \sym (X , \wa \rho_\bullet)$; since $\autrb = \left\{ id \right\}$, from the previous lemma we conclude $\wa \rho \simeq X \wa \rho_\bullet$. We denote by $\mP_{X,d} \simeq X \wa \ud$ the unique (up to isomorphism) element of $\sym ( X , \wa \ud )$.
\end{ex}

\subsubsection{Special categories.}
Let $G \subseteq \ud$ be a closed group. We denote by $NG$ the normalizer of $G$ in $\ud$, and by $QG := NG/G$ the quotient group, so that we have an epimorphism 
\begin{equation}
\label{def_p}
p : NG \to QG \ .
\end{equation}
\noindent The inclusion $G \subseteq \ud$ implies that $\mO_\ud \subseteq \mO_G$. We denote by ${\bf aut} ( \mO_d , \mO_G )$ $\subseteq$ ${\bf aut}_\theta \mO_d$ the group of automorphisms of $\mO_d$ leaving $\mO_G$ globally stable, and coinciding with the identity on $\mO_\ud$. From \cite[Cor.3.3]{DR87}, we conclude that ${\bf aut} ( \mO_d , \mO_G )$ is isomorphic to a subgroup of $\ud$, acting on $\mO_d$ according to (\ref{def_cag}). 

\begin{thm}
\label{thm_qg}
Let $( \rho , d , V )$ be a special object, and suppose $(\ii) = \bC$. Then, there is a closed group $G \subseteq \sud$ with a group isomorphism ${\bf aut}_\eps \wa \rho \simeq QG$. Let $\pi : {\bf aut} ( \mO_d , \mO_G  ) \to  {\bf aut} \mO_G$ be the map assigning to $\alpha \in  {\bf aut} ( \mO_d , \mO_G )$ the restriction on $\mO_G$; then, there is a commutative diagram of group morphisms
\begin{equation}
\label{cd_qg}
\xymatrix{
   NG
   \ar[r]^-{ p }
   \ar[d]_-{ \wa{} \ }
 & QG
   \ar[d]^-{ \ \wa{} }
\\  {\bf aut} ( \mO_d , \mO_G )
   \ar[r]^-{ \pi  }
 & {\bf aut}_\theta \mO_G
}
\end{equation}
\noindent The vertical arrows of the above diagram are group isomorphisms.
\end{thm}

\begin{proof}
By Lemma \ref{thm_spec}, there are isomorphisms ${\bf aut}_\eps \wa \rho \simeq {\bf aut}_\theta \wa G \simeq {\bf aut}_\theta \mO_G$; moreover, for every $\alpha \in {\bf aut}_\eps \wa \rho$ there exists $u \in \ud$ such that $\wa u \circ i = i \circ \alpha$. Now $i (\oro) = \mO_G$, so that the previous equality implies that $\wa u \in {\bf aut} \mO_d$ restricts to an automorphism of $\mO_G$; moreover, for every $g \in G$ we find
\begin{equation}
\label{eq_qg}
\wa u \circ \wa g \circ \wa{u^*} (t) = 
\wa u \circ \wa g (\wa{u^*}(t)) =
\wa u \circ \wa{u^*} (t) = t
\end{equation}
\noindent (we used $\wa{u^*} (t) \in \mO_G$, and $\wa g \in {\bf aut}_{\mO_G} \mO_d$). Thus, we conclude that $\wa u \circ \wa g \circ \wa{u^*} \in {\bf aut}_{\mO_G} \mO_d$, i.e. $ugu^* = g'$ for some $g' \in G$; in other terms, $u \in NG$. Moreover, it is clear that $\wa{ug} \circ i = \wa u \circ i = i \circ \alpha$, thus $u \in NG$ is defined up to multiplication by elements of $G$. Moreover, note that $\wa u |_{\mO_G} = \wa v |_{\mO_G}$, $u,v \in NG$, implies $\wa{uv^*} \in {\bf aut}_{\mO_G} \mO_d$, i.e. $uv^* \in G$. On the converse, if $u \in NG$, $t \in \mO_G$, then $\wa g \circ \wa u (t) = \wa u \circ \wa{g'} (t) = \wa u (t)$ for some $g' \in G$, and this implies $\wa u \in {\bf aut} (\mO_d , \mO_G)$. We conclude that the map $\pi (\wa u) \mapsto p(u)$, $\wa u \in {\bf aut} ( \mO_d , \mO_G )$, is an isomorphism. 
\end{proof}

\begin{defn}
\label{def_spec_cat}
A {\bf special category} is a locally trivial, symmetric tensor \sC category $( \wa \rho , \otimes , \iota , \eps )$ with fibre $( \wa \rho_\bullet , \otimes_\bullet , \iota_\bullet , \eps_\bullet )$, such that $\rho_\bullet$ is a special object.
\end{defn}

We emphasize the fact that a special category is locally trivial in the sense of Def.\ref{def_slc}, thus the local charts preserve tensor product and symmetry. This also implies that $(\iota_\bullet , \iota_\bullet) \simeq \bC$.

\begin{thm}
\label{thm_amen}
Let $( \wa \rho  , \otimes , \iota , \eps )$ be a special category with fibre $( \wa \rho_\bullet , \otimes_\bullet , \iota_\bullet , \eps_\bullet )$. Then, $\rho$ is amenable. There exists $d \in \bN$ and a unique up to isomorphism compact Lie group $G \subseteq \sud$ such that $\wa \rho_\bullet \simeq \wa G$. Thus, $\oro$ is a locally trivial bundle of \sC algebras with fibre $\mO_G$, and there is a one-to-one correspondence 
\begin{equation}
\label{eq_h1}
\sym ( X^\iota , \wa \rho_\bullet ) \simeq H^1(X^\iota , QG) \ .
\end{equation}
\end{thm}

\begin{proof}
It follows from point (3) of Lemma \ref{thm_spec} that the fibre of $\wa \rho$ is isomorphic to $(  \wa G , \otimes  , \iota_\bullet , \theta  )$, where $G \subseteq \sud$ is a compact Lie group unique up to isomorphism. Note that $\wa G$ is amenable (see Ex.\ref{ss_cog}). Let now $\rho_* \in {\bf end} \oro$ be the canonical endomorphism, so that $(\rhors) \subseteq ( \rho_*^r , \rho_*^s )$, $r,s \in \bN$. We consider a cover of closed neighborhoods $\left\{  X_i \right\}$ with local charts $\pi_i : \wa \rho \to X_i \wa G$, and a partition of unity $\left\{ \lambda_i \right\}$ subordinate to $\left\{ \right. \dot{X}_i \left. \right\}$. Since $\wa G$ is amenable, every $X_i \wa G$ is amenable; thus, if $t \in ( \rho_*^r , \rho_*^s )$, then $\pi_i ( \lambda_i t) \in  C_0 ( X_i , (\sigma_G^r , \sigma_G^s) ) = C_0( X_i , (\hrs)_G )$. Since $\pi_i^{-1} ( C_0( X_i , (\hrs)_G ) )$ is contained in $(\rhors)$, we conclude that $\lambda_i t \in (\rhors)$. Thus $t = \sum_i \lambda_i t \in (\rhors)$, and $\rho$ is amenable. By Thm.\ref{thm_qg} and Lemma \ref{lem_lc} we obtain (\ref{eq_h1}). The assertions about $\oro$ follow from Prop.\ref{main_thm}. 
\end{proof}

For every $\mT \in$ $\sym ( X , \wa G )$, we denote by $Q \mT \in$ $H^1 (X,QG)$ the unique-up-to-isomorphism principal $QG$-bundle associated with $\mT$. The map (\ref{eq_h1}) has to be intended as an isomorphism between sets endowed with a distinguished element, in the sense that $Q ( X \wa G )$ coincides with the trivial principal $QG$-bundle $X \times QG$. 

Some particular cases follow.
If $QG$ is Abelian, then $H^1 ( X , QG )$ is a sheaf cohomology group (\cite[I.3.1]{Hir}); this allows one to define a group structure on $\sym ( X , \wa G )$.
Let $SY$ be the suspension of a compact Hausdorff space $Y$, and $QG$ arcwise connected. By classical arguments (\cite[Chp.7.8]{Hus}), we have an isomorphism $\sym ( SY , \wa G ) \simeq [ Y , QG ]$, where $[ Y , QG ]$ is the set of homotopy classes of continuous maps from $Y$ into $QG$. In particular, if $Y = S^n$ is the $n$-sphere (i.e. $X = S^{n+1}$), then $\sym ( S^{n+1} , \wa G )$ is isomorphic to the homotopy group $\pi_n ( QG )$.

\begin{prop}
\label{prop_dc}
Let $( \wa \rho , \otimes , \iota , \eps )$ be a special category. Then, there are $d \in \bN$ and $\mV \subset ( \iota , \rho^d )$ such that the triple $( \rho , d , \mV )$ defines a twisted special object.
\end{prop}

\begin{proof}
We take $d$ as the dimension of the fibre $\rho_\bullet$ in the sense of (\ref{def_dim}). Let us consider the totally antisymmetric projection $P_{\rho , \eps , d} \in (\rho^d , \rho^d)$. To be concise, we also write $P_{\theta,d} := P_{\bH_d , \theta , d} \in ( \bH_d^d , \bH_d^d )$; if $\pi_i : \wa \rho \to X_i \wa G$ is a local chart in the sense of (\ref{def_lc}), then $\pi_i ( P_{\rho, \eps , d} ) = P_{\theta , d}$, in fact $\pi_i (\eps_\rho (p)) = \theta(p)$, $p \in \bP_d$. Now, $P_{\theta , d} \in ( \bH_d^d , \bH_d^d )_G$ and has rank one, in fact $P_{\theta , d} = S \circ S^*$, where $S \in ( \iota_\bullet , \bH_d^d )_G$ satisfies $(S^* \otimes 1_d ) \circ (1_d \otimes S) = (-1)^{d-1} d^{-1} 1_d$ (\cite[Lemma 2.2]{DR87}). We define
\[
\mV :=
\left\{ 
V \in ( \iota , \rho^d ) : V = P_{\rho , \eps , d} \circ V
\right\} \ .
\]
\noindent It is clear that $\mV$ is a vector space. Moreover, $\mV$ has a natural structure of a Hilbert $C(X^\iota)$-bimodule, by considering the actions $f , V \mapsto$ $V \circ f$, $f \otimes V$, $f \in C(X^\iota) = (\ii)$. We now consider a partition of unity $\left\{ \lambda_i \right\} \subseteq (\ii)$, subordinate to the open cover $\left\{ \dot{X}_i \right\}$; if $V \in \mV$, then 
\[
\pi_i ( \lambda_i^{1/2} V) = 
\pi_i ( \lambda_i^{1/2} P_{\rho, \eps , d} \circ V ) =
P_{\theta , d} \circ \pi_i ( \lambda_i^{1/2} V ) \ .
\]
\noindent Thus, there exists $f_i \in C(X^\iota)$ such that $\pi_i ( \lambda_i^{1/2} V ) = f_i S$; if $V' \in \mV$, with $\pi_i ( \lambda_i^{1/2} V' ) = f'_i S$, then $\lambda_i V^* \circ V' = f_i^* f'_i$. This implies
\[
\begin{array}{ll}
\lambda_i ( V^* \otimes 1_\rho ) \circ ( 1_\rho \otimes V' ) & =
\pi_i^{-1} 
(   f_i^* f'_i  \ 
( S^* \otimes_\bullet 1_d ) 
\circ ( 1_d \otimes_\bullet S  
) ) \\ & =
(-1)^{d-1} d^{-1} \ f_i^* f'_i \ 1_d     \\ & =
\lambda_i (-1)^{d-1} d^{-1} ( V^* \circ V' ) \otimes 1_\rho \ .
\end{array}
\]
\noindent By summing over the index $i$, we obtain (\ref{eq_tso}(2)). In the same way, the equality $\pi_i ( \lambda_i V' \circ V^*) = f'_i f_i^* P_{\theta , d}$ implies (\ref{eq_tso}(3)). 
\end{proof}

Let $u \in NG$. Since $G \subseteq \sud$, we find $\det (ug) = \det (u)$ for every $g \in G$. This means that the determinant factorizes through a morphism $\det_Q : QG \to \bT$. By functoriality of $H^1(X,\cdot \ )$, a map $\det_{Q,*} : H^1 (X,QG)$ $\to$ $H^1(X,\bT) \simeq$ $H^2(X,\bZ)$ is induced.

\begin{cor}
\label{cor_ch_sc}
Let $( \wa \rho , \otimes , \iota , \eps )$ be a special category with associated $QG$-cocycle $Q \wa \rho \in$ $H^1(X,QG)$. Then, $\rho$ has Chern class $c(\rho) =$ $\det_{Q,*} (Q \wa \rho)$.
\end{cor}

\begin{proof}
For every $y \in QG$, we denote by $\wa y \in$ ${\bf aut}_\theta \wa G \simeq$ ${\bf aut}_\theta \mO_G$ the associated autofunctor defined according to (\ref{cd_qg}). We adopt the same notation if $y : X \to QG$ is a continuous map (so that, $\wa y$ is a continuous ${\bf aut}_\theta \wa G$-valued map).
Let $Q :=$ $( \left\{ X_i \right\} , \left\{ y_{ij} \right\})$ be a $QG$-cocycle associated with $\wa \rho$. If $\pi_i : \wa \rho \to X_i \wa G$ are local charts associated with $Q$ (i.e., $\pi_i \circ \pi_j^{-1} =$ $\wa y_{ij}$), then it follows from the proof of Prop.\ref{prop_dc} that $\pi_i (\mV) =$ $C(X_i) \otimes \wedge^d \bH_d$, where $\wedge^d \bH_d := P_{\theta,d} H_d^d$ is the rank one Hilbert space of totally antisymmetric vectors. Now, if $u \in \ud$ then $\wa u (S) =$ $\det u \cdot S$ for every $S \in \wedge^d \bH_d$; so that, $\wa y (S) =$ $\det_Q (y) \cdot S$, $y \in QG$. Thus, for every $V \in$ $C(X_{ij}) \otimes \wedge^d \bH_d$ we find
\[
\pi_i \circ \pi_j^{-1} (V) =
\wa y_{ij} (V) =
{\det}_Q (y_{ij}) \cdot V \ .
\]
\noindent This implies that the line bundle associated with $\mV$ has transition maps $\det_Q (y_{ij}) :$ $X_{ij} \to $ $\bT$. 
\end{proof}

\begin{ex}
\label{ex_sud}
Let us consider the case $G = \sud$, so that $QG = \bT$ and (\ref{def_p}) is the determinant map $\det : \ud \to \bT$. By Thm.\ref{thm_amen}, we obtain 
\[
\sym (X,\wa \sud) = H^1(X,\bT) \simeq H^2(X,\bZ) \ ,
\]
\noindent and the map $\det_{Q,*}$ defined in Cor.\ref{cor_ch_sc} is the identity. 
The previous considerations imply that elements of $\sym (X,\wa \sud)$ are labeled by pairs $(d,\zeta)$, $d \in \bN$, $\zeta \in H^2(X,\bZ)$. We denote by $\mT_{d,\zeta}$ the generic element of $\sym (X,\wa \sud)$, and by $( \rho , d , \mV )$ the twisted special object generating $\mT_{d,\zeta}$. By Cor.\ref{cor_ch_sc}, the class $\zeta$ is the first Chern class of the line bundle $\mL (\zeta) \to X$ associated with $\mV$. 
The arrows of $\mT_{d,\zeta}$ are generated by the symmetries $\eps_\rho (r,s)$, $r,s \in \bN$, and elements of $\mV$, in the same way as the dual $\wa{\sud}$ is generated by the flip operators $\theta (r,s)$ and the totally antisymmetric isometry $S \in$ $(\iota,\bH_d^d)$ (see Ex.\ref{ss_cog}, \cite[Lemma 3.7]{DR87})
\end{ex}

\section{A classification for certain $\mO_G$-bundles.}
\label{ss_ogb}

We recall the reader to the notation of Ex.\ref{ss_cog}. 
Let $G \subseteq \sud$ be a compact group, and $\wa \rho \in \sym ( X, \wa G )$. Since $\rho$ is amenable (Lemma \ref{lem_oro_tso}, Prop.\ref{prop_dc}) we find $\wa \rho_* = \wa \rho$, i.e. $( \rhors ) = ( \rho_*^r , \rho_*^s  )$, $r,s \in \bN$. By Prop.\ref{main_thm}, the DR-algebra $\oro$ is a locally trivial continuous bundle with fibre $\mO_G$ (an $\mO_G$-{\em bundle}, for brevity). The structure group of $\oro$ is clearly given by ${\bf aut}_\theta \mO_G \simeq QG$, in the sense that $\oro$ admits a set of transition maps taking values in ${\bf aut}_\theta \mO_G \subseteq$ ${\bf aut} \mO_G$. 

On the converse, let $\mA$ be an $\mO_G$-bundle with structure group ${\bf aut}_\theta \mO_G$. We consider an ${\bf aut}_\theta \mO_G$-cocycle $( \left\{ X_h \right\}_h^n , \left\{ \alpha_{hk} \right\} )$ associated with $\mA$, such that each $X_h$ is a closed neighborhood. We denote by $\alpha_{hk}^x \in$ ${\bf aut}_\theta \mO_G$ the evaluation of $\alpha_{hk}$ over $x \in X_{hk}$. By construction, there is a one-to-one correspondence between elements $a \in \mA$ and $n$-ples
\begin{equation}
\label{eq_ccA}
\begin{array}{ll}
( a_h  ) \in \bigoplus_h^n  \left( C(X_h) \otimes \mO_G \right) \ \ : \ \
a_h (x) = \alpha_{hk}^x ( a_k (x) )  \ , \ x \in X_{hk}
\ .
\end{array}
\end{equation}
\noindent Now, by definition of ${\bf aut}_\theta \mO_G$ we have  $\sigma_G \circ \alpha_{hk}^x =$ $\alpha_{hk}^x \circ \sigma_G$, $x \in X_{hk}$. Let us denote by $\iota_h$ the identity automorphism on $C(X_h)$; by the above considerations, for every $a \in \mA$ the $n$-ple $( ( \iota_h \otimes \sigma_G ) (a_h) )$ satisfies (\ref{eq_ccA}), and defines an element of $\mA$, denoted by $\rho (a)$. An immediate check shows that the map $\left\{ a \mapsto \rho (a) \right\}$ is a $C(X)$-endomorphism of $\mA$; we call $\rho$ the {\em canonical endomorphism} of $\mA$. By construction, we have a set of local charts 
\[
\pi_h : ( \mA , \rho ) \to ( \ C(X_h) \otimes \mO_G \ , \  \iota_h \otimes \sigma_G \ )
\ , \ \ 
\pi_h (a) := a_h \ ,
\]
\noindent so that in particular 
\begin{equation}
\label{eq_phrs}
\pi_h (\rhors) =
C(X_h) \otimes (\sigma_G^r , \sigma_G^s)
\ , \ \
r,s \in \bN \ .
\end{equation}
\noindent Since $\wa \sigma_G = \wa G$, the tensor category $\wa \rho \subseteq$ ${\bf end}_X \mA$ (defined as in Ex.\ref{ex_end}) is an element of $\sym ( X, \wa G  )$. Moreover, $\mA$ is the DR-algebra associated with $\rho$, in fact the set $\left\{ (\rhors) \right\}_{r,s}$ is total in $\mA$ (this is easily verified by using (\ref{eq_phrs}), and the fact that $\left\{ (\sigma_G^r , \sigma_G^s) \right\}_{r,s}$ is total in $\mO_G$). Thus, there is a unique (up-to-equivalence) $QG$-cocycle associated with $(\mA,\rho)$, namely
\[
Q (\mA,\rho) := Q \wa \rho \ \in \  H^1 ( X,QG ) \ .
\]
\noindent Let $\mA$, $\mA'$ be $\mO_G$-bundles  with structure group ${\bf aut}_\theta \mO_G$. We denote by $\rho \in {\bf end}_X \mA$, $\rho' \in {\bf end}_X \mA'$ the canonical endomorphisms, and by $\eps \in (\rho^2 , \rho^2)$, $\eps' \in ( {\rho'}^2 , {\rho'}^2 )$ the flip operators defined by the associated symmetries. Moreover, $\otimes$ (resp. $\otimes'$) denotes the tensor product in ${\bf end}_X \mA$ (resp. ${\bf end}_X \mA'$), and $\iota \in {\bf aut}_X \mA$ (resp. $\iota' \in {\bf aut}_X \mA'$) the identity automorphism. By recalling (\ref{eq_aut}), we obtain a translation of Thm.\ref{thm_amen} in terms of \sC algebra bundles:

\begin{prop}
\label{prop_cog}
With the above notation, the following are equivalent:
\begin{enumerate}
\item  $Q (\mA,\rho) = Q (\mA',\rho')$ $\in$ $H^1 ( X,QG )$;
\item  there is a $C(X)$-isomorphism $\alpha : ( \mA , \rho ) \to$ $( \mA' , \rho' )$, with $\alpha (\eps) = \eps'$;
\item there is an isomorphism $\alpha_* :$ $( \wa \rho , \otimes , \eps , \iota ) \to$ $( \wa \rho' , \otimes' , \eps' , \iota' )$.
\end{enumerate}
\end{prop}

We recall the reader to Ex.\ref{ex_sud}, and briefly discuss the case $G = \sud$. We consider $\mT_{d,\zeta} \in$ $\sym ( X , \wa{\sud} )$, $d \in \bN$, $\zeta \in$ $H^2(X,\bZ)$; moreover, we denote by $\mO_{d,\zeta}$ the DR-algebra associated with the generating twisted special object $( \rho , d , \mV )$. Now, Ex.\ref{ex_pc} implies that the \sC subalgebra $\mO_\eps \subset$ $\mO_{d,\zeta}$ generated by $\left\{ \eps_\rho (r,s) \right\}_{r,s}$ is isomorphic to $C(X) \otimes \mO_\ud$. Let us consider the (graded) endomorphism $\nu \in$ ${\bf end} \mO_{d,\zeta}$ defined as in (\ref{def_nu}); since $\mO_{d,\zeta}^0 = \mO_\eps$, we find that $\nu$ restricts to an endomorphism of $\mO_\eps$, that we regard as an endomorphism $\phi \in $ ${\bf end} \left( C(X) \otimes \mO_\ud \right)$. By Ex.\ref{ex_sud} it follows that $\mV$ is the module of sections of the line bundle $\mL (\zeta) \to X$ with first Chern class $\zeta$, and that $\mO_{d,\zeta}$ is generated by $\mO_\eps$ and $\mV$; thus, (\ref{def_nu}) implies that $\mO_{d,\zeta}$ is a crossed product
\[
\mO_{d,\zeta} \simeq \left( C(X) \otimes \mO_\ud \right) \rtimes_\phi^{\mL(\zeta)} \bN \ ,
\]
\noindent in the sense of \cite[\S 3]{Vas05}.
In the case in which $X$ is the 2-sphere $S^2$, it is well-known that $H^2(S^2,\bZ) \simeq \bZ$. By Prop.\ref{prop_cog}, the set of isomorphism classes of $\mO_\sud$-bundles with structure group ${\bf aut}_\theta \mO_\sud$ and base space $S^2$ is labeled by $\bZ$; in particular, $0 \in \bZ$ corresponds to the trivial bundle $C(S^2) \otimes \mO_\sud$.

\section{Outlooks.}

A duality theory for special categories will be the next step w.r.t. the present work. In explicit terms, our aim is to prove a generalization of Lemma \ref{thm_spec} (i.e., \cite[Thm.4.17]{DR89}): instead of the dual $\wa G$, our model category will be the one of tensor powers of a vector bundle $\mE \to X$, with arrows morphisms equivariant w.r.t. a group bundle $\mG \to X$. The $\mG$-action on the tensor powers $\mE^r$, $r \in \bN$, will be defined according to the gauge-equivariant $K$-theory introduced in \cite[\S 1]{NT04}, so that $\mG$ will play the role of a dual object. 

An important step to prove such a duality is to find an embedding functor in the sense of Def.\ref{def_ef}. As outlined in the introduction of the present paper, existence and unicity of the embedding functor (and the dual object $\mG$) are not ensured, in contrast with the case $(\ii) \simeq \bC$. In particular, non-isomorphic dual objects may be associated with the same special category, in the case in which more than an embedding functor exists.

We will use the cohomological classification (\ref{eq_h1}) to provide a complete description of such a phenomenon in geometrical terms. Given a special category $\mT \in$ $\sym ( X,\wa G )$ and the associated principal $QG$-bundle $Q \mT \in$ $H^1(X,QG)$, our tasks will be the following:
\begin{enumerate}
\item determine which are the geometrical properties required for $Q \mT$ in order to find an embedding functor $i : \mT \hra {\bf vect}(X)$;
\item in the case in which there exists the embedding $i$, give a characterization of $\mE$ and the dual object $\mG$ in terms of geometrical properties of $Q \mT$, and determine the dependence of $\mE$, $\mG$ on $i$.
\end{enumerate}

These questions have a natural translation in terms of \sC algebra bundles and \sC dynamical systems. If $(\mA,\rho)$ is a \sC dynamical system as in Sec.\ref{ss_ogb}, one could ask whether there exists a vector bundle $\mE \to X$ with associated DR-dynamical system $( \coe , \bT , \sigma )$ (see Ex.\ref{ex_od}), such that there is a $C(X)$ - monomorphism $\phi : (\mA,\rho) \hra$ $(\coe,\sigma)$. The existence of an embedding functor for $\wa \rho$ is equivalent to the existence of $\phi$, and the obstrution for the existence can be encoded by a cohomological invariant.

\noindent {\bf Acknowledgements.} The author wishes to thank an anonymous referee for many precious remarks on a previous version of the present work, in particular for the references \cite{Kan01,DG,BK04}.



\begin{thebibliography}{99}

%
%

\bibitem{APT73}
C.A. Akemann, G.K. Pedersen, J. Tomiyama, Multipliers in \sC algebras. J. Funct. Anal. 13 (1973) 277-301.

\bibitem{Ati}
M.F. Atiyah, $K$-Theory, Benjamin, New York (1967).

%
%


\bibitem{Bla96}
E. Blanchard, D\'eformations de \sC alg\'ebres de Hopf, {Bull. Soc. math. France} {124} (1996) 141-215.

\bibitem{BK04}
E. Blanchard, E. Kirchberg. Global Glimm halving for {\it C*}-bundles. J.
Oper. Theory 52 (2004) 385-420.


\bibitem{BL04}
H. Baumg\"artel, F. Lled\'o, Duality of compact groups and Hilbert {\it C*}-systems for {\it C*}-algebras with a nontrivial center, Int. J. Math. 15 (8) (2004) 759-812.



\bibitem{Cun77}
J. Cuntz, Simple {\it C*}-algebras Generated by Isometries, Comm.
Math. Phys. 57 (1977) 173-185.


\bibitem{Del90}
P. Deligne, Categories tannakiennes, in The Grothendieck Festschrift Volume II,
Cartier, P., et al., Birkhauser, Boston, (1990) 111-196.

\bibitem{DG}
M.J. Dupr\'e, R.M. Gillette. Banach bundles, Banach modules and automorphisms
of {\it C*}-algebras. Pitman (1983).


\bibitem{Dix}
J. Dixmier: {\it C*}-algebras, North-Holland Publishing Company, Amsterdam - New
York, Oxford (1977).


\bibitem{DPR01}
S. Doplicher, C. Pinzari, J.E. Roberts, An Algebraic Duality Theory for Multiplicative Unitaries, Int. J. Math. 12 415-459 (2001).


\bibitem{DR87}
S. Doplicher, J.E. Roberts, Duals of Compact Lie Groups Realized in the Cuntz Algebras and Their Actions on {\it C*}-algebras, {J. Funct. Anal.} {74} (1987) 96-120.


\bibitem{DR89}
S. Doplicher, J.E. Roberts, A New Duality Theory for Compact Groups, {Invent. Math.} {98} (1989) 157-218.



\bibitem{DR90}
S. Doplicher, J.E. Roberts, Why there is a field algebra with a compact gauge group describing the superselection structure in particle physics, Commun. Math. Phys. 131 (1990) 51-107.


\bibitem{DR89A}
S. Doplicher, J.E. Roberts, Endomorphisms of \sC algebras, Cross Products and Duality for Compact Groups, {Ann. Math.} {130} (1989) 75-119.



\bibitem{Dup74}
M.J. Dupr\'e, Classifying Hilbert bundles I, J. Funct. Anal. 15 (1974) 244-278.


%
%


\bibitem{Hir}
F. Hirzebruch, Topological Methods in Algebraic Geometry, Springer-Verlag, 1966.


\bibitem{Hus}
D. Husemoller: Fiber Bundles, Mc Graw-Hill Series in Mathematics, 1966.


\bibitem{KPW04}
T. Kajiwara, C. Pinzari, Y. Watatani, Jones index theory for Hilbert {\it C*}-bimodules and its equivalence with conjugation theory, J. Funct. Anal. 215 (1) (2004) 1-49.


\bibitem{Kan01}
T. Kandelaki. Multiplier and Hilbert {\it C*}-categories. Georgian Academy of
Sciences. Proceedings of A. Razmadze Mathematical Institute 127 (2001) 89-111.


\bibitem{Kar}
M. Karoubi: {$K$-Theory}, Springer Verlag, Berlin - Heidelberg - New York, 1978.


\bibitem{Kas88}
G.G. Kasparov, Equivariant $KK$-Theory and the Novikov Conjecture, {Invent. Math.} {91} (1988) 147-201.


\bibitem{KW95}
E. Kirchberg, S. Wassermann, Operations on Continuous Bundles of \sC algebras, Math. Annalen 303 (1995) 677-697.



\bibitem{LR97}
R. Longo, J.E. Roberts, A Theory of Dimension, {$K$-Theory} {11} (1997) 103-159.

%
%

\bibitem{Mit01}
P.D. Mitchener, Symmetric $K$-theory Spectra of {\it C*}-categories, {$K$-theory} {24}(2) (2001) 157-201.


\bibitem{NT04}
V. Nistor, E. Troitsky, An index for gauge-invariant operators and the Dixmier-
Douady invariant, Trans. AMS. 356 (2004) 185-218


\bibitem{Nil96}
M. Nilsen, { {\it C*}-Bundles and $C_0(X)$-algebras}, Indiana Univ. Math. J. {45} (1996) 463-477.


\bibitem{Ped0}
G.K. Pedersen, Analysis Now. Springer Verlag (1989).


\bibitem{Ped}
G.K. Pedersen, {\it C*}-algebras and their automorphism groups, London Math. Soc.
Monographs, Academic Press, London-New York (1979).


%
%


\bibitem{Seg68}
G. Segal, Equivariant $K$-theory, Inst. Hautes \'Etudes Sci. Publ. Math. 34 (1968) 129-151.

%
%

\bibitem{Vas}
E. Vasselli, Continuous fields of {\it C*}-algebras Arising from Extensions of Tensor {\it C*}-categories, {J. Funct. Anal.} {199} (2003) 122-152.


\bibitem{Vas05}
E. Vasselli, Crossed Products by Endomorphisms, Vector Bundles and Group Duality, Int. J. Math. 16 (2) (2005) 137-171.

\bibitem{Vas05A}
E. Vasselli, Crossed Products by Endomorphisms, Vector Bundles and Group Duality, II, Int. J. Math. 17(1) (2006) 65-96.


\bibitem{Vas05B}
E. Vasselli, The {\it C*}-algebra of a vector bundle and fields of Cuntz algebras, J. Funct. Anal. 222(2) (2005) 491-502.


\bibitem{Vas06}
E. Vasselli, Group bundle duality, invariants for certain \sC algebras, and twisted equivariant $K$-theory, preprint arXiv math.KT/0605114 (2006), to appear on Proceedings of the Conference {\it C*}-algebras and elliptic theory, II.


\bibitem{Wor87}
S.L. Woronowicz, Tannaka-Krein Duality for Compact Matrix Pseudogroups. Twisted $SU(N)$-groups, {Invent. Math.} {111} (1987) 35-76.


\bibitem{Zit05}
P. Zito, 2-\sC categories with non-simple units, Preprint arXiv math.CT/0509266 (2005), to appear on Adv. Math..


\end{thebibliography}
\end{document}